\documentclass[12pt]{amsart}
\usepackage{pstricks}
\usepackage{latexsym}
\usepackage{amssymb}

\numberwithin{section}{part}

 \textwidth=\textwidth \textheight=\textheight
\theoremstyle{plain}
\newtheorem{thm}{Theorem}
\newtheorem{lem}[thm]{Lemma}
\newtheorem{prop}[thm]{Proposition}
\newtheorem{corollary}[thm]{Corollary}

\newtheorem{observation}[thm]{Observation}
\theoremstyle{definition}
\newtheorem{definition}[thm]{Definition}

\newtheorem{remark}[thm]{Remark}

\newtheorem{example}[thm]{Example}

\subjclass[2000]{05E, 17B}

\date{November~9, 2003}

\newcommand{\set}[1]{ \left\{ #1 \right\} }
\newcommand{\ov}[1]{\overline{#1}}
\def\P{{\mathcal P}}
\def\u{{\mathcal U}}
\def\tu{\tilde{{\mathcal U}}}
\def\v{{\mathcal V}}
\def\tv{\tilde{{\mathcal V}}}
\def\a{\alpha}
\def\b{\beta}

\def\ww{{\mathbf w}}
\def\C{\mathbb C}
\def\ll{\lambda}
\def\tll{\tilde{\lambda}}
\def\d{\delta}

\def\p{{\mathcal P}}

\def\zz{{\mathbb Z}}
\def\S{{\mathcal S}}
\def\Chi{{\mathbf \chi}}
\def\uqsln{U_q(\widehat{ {\mathfrak
sl}}_n)}

\def\g{{\mathcal G}}
\def\f{{\mathbf F}}

\def\CH{{\C(q)[H_-]}}
\def\h{{\mathcal H}}
\def\hh{{\mathbf h}}
\def\e{{\mathbf e}}

\def\m{{\mathbf \Phi}}
\def\o{{\Omega}}
\def\to{\tilde{\Omega}}
\def\l{{\mathcal L}}
\def\xx{{\mathcal X}}
\def\kk{{\mathcal K}}
\def\plqn{{\Upsilon_{q,n}}}
\def\w{{\omega}}
\def\log{{\rm log}}

\newcommand{\diff}[1]{\frac{\partial}{\partial #1}}
\newcommand{\ip}[1]{\left\langle #1 \right \rangle}
\newcommand{\be}[1]{\begin{equation} \label{#1}}
\newcommand{\ee}{\end{equation}}

\newcommand{\br}[1]{|#1\rangle}
\newcommand{\brac}[1]{\left(#1\right)}

\psset{unit=1pt} \psset{arrowsize=4pt 1} \psset{linewidth=1pt}
\newrgbcolor{mygreen}{.2 .7 .2}
\newrgbcolor{myred}{.7 .3 .2}
\newgray{mygray}{.90}
\countdef\x=23 \countdef\y=24 \countdef\z=25 \countdef\t=26

\def\abox(#1,#2)#3{
\x=#1 \y=#2 \multiply\x by 16 \multiply\y by 16 \z=\x \t=\y
\advance\z by 16 \advance\t by 16
\psline(\x,\y)(\x,\t)(\z,\t)(\z,\y)(\x,\y) \advance\x by 8
\advance\y by 8 \rput(\x,\y){{\bf #3}}}

\def\hdom(#1,#2)#3{
\x=#1 \y=#2 \multiply\x by 16 \multiply\y by 16 \z=\x \t=\y
\advance\z by 32 \advance\t by 16
\psline(\x,\y)(\x,\t)(\z,\t)(\z,\y)(\x,\y) \advance\x by 16
\advance\y by 8 \rput(\x,\y){{\bf #3}}}

\def\vdom(#1,#2)#3{
\x=#1 \y=#2 \multiply\x by 16 \multiply\y by 16 \z=\x \t=\y
\advance\z by 16 \advance\t by 32
\psline(\x,\y)(\x,\t)(\z,\t)(\z,\y)(\x,\y) \advance\x by 8
\advance\y by 16 \rput(\x,\y){{\bf #3}}}

\def\rec(#1,#2,#3,#4){
\psline(#1,#2)(#3,#2)(#3,#4)(#1,#4)(#1,#2) }

\begin{document}
\title{Ribbon Tableaux and the Heisenberg Algebra}
\author{Thomas Lam}
\address{Department of Mathematics,
         M.I.T., Cambridge, MA 02139}

\email{thomasl@math.mit.edu}
\begin{abstract}
In \cite{LLT} Lascoux, Leclerc and Thibon introduced symmetric
functions $\g_\ll$ which are spin and weight generating
functions for ribbon tableaux.  This article is aimed at studying
these functions in analogy with Schur functions.  In particular we
will describe:
\begin{itemize}
\item
a Pieri and dual-Pieri formula for ribbon functions,
\item
a ribbon Murnagham-Nakayama formula,
\item
ribbon Cauchy and dual Cauchy identities,
\item
and a $\C$-algebra isomorphism $\w_n: \Lambda(q) \rightarrow
\Lambda(q)$ which sends each $\g_\ll$ to $\g_{\ll'}$.
\end{itemize}
Our study of the functions $\g_\ll$ will be connected to the Fock
space representation $\f$ of $\uqsln$ via a linear map $\m: \f
\rightarrow \Lambda(q)$ which sends the standard basis of $\f$ to
the ribbon functions.  Kashiwara, Miwa and Stern \cite{KMS} have
shown that a copy of the Heisenberg algebra $H$ acts on $\f$
commuting with the action of $\uqsln$.  Identifying the Fock Space
of $H$ with the ring of symmetric functions $\Lambda(q)$ we will
show that $\m$ is in fact a map of $H$-modules with remarkable
properties. We give a combinatorial proof that the ribbon
Murnagham-Nakayama and Pieri rules are formally equivalent thus
allowing us to describe the action of the genrators of $H$ on $\f$
in terms of `border ribbon strips'. We will also connect the
ribbon Cauchy and Pieri formulae to the combinatorics of ribbon
insertion as studied by Shimozono and White \cite{SW2}.  In
particular we give complete combinatorial proofs for the domino
$n=2$ case.
\end{abstract}
\maketitle \tableofcontents

\part*{Introduction}
Let $n \geq 1$ be a fixed integer and $\ll$ a partition with empty $n$-core.
In analogy with the combinatorial
definition of the Schur functions, Lascoux, Leclerc and Thibon \cite{LLT}
have defined a family of symmetric functions $\g_\ll(X;q) \in \Lambda(q)$ by:
\[
\g_\ll(X;q) = \sum_T q^{s(T)}{\mathbf x}^{w(T)}
\]
where the sum is over all \emph{semistandard ribbon tableaux} of shape $\ll$,
and $s(T)$ and $w(T)$ are the spin and weight of $T$ respectively.
The definition of a semistandard ribbon tableaux is analagous to the definition
of semistandard Young tableaux, with boxes replaced by ribbons (or border strips)
of length $n$.  We shall loosely call the functions $\g_\ll(X;q)$
\emph{ribbon functions}.

\begin{figure}[ht]
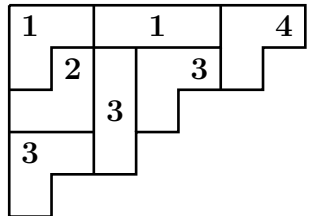

\pspicture(0,20)(120,80)

\psline(0,90)(112,90)(112,74)(96,74)(96,58)(64,58)(64,42)(48,42)(48,26)(16,26)(16,10)(0,10)(0,90)
\psline(0,58)(16,58)(16,74)(80,74)(80,90) \psline(80,74)(80,58)
\psline(32,90)(32,26) \psline(48,74)(48,42) \psline(0,42)(32,42)
\rput(8,34){{\bf 3}} \rput(24,66){{\bf 2}} \rput(40,50){{\bf 3}}
\rput(72,66){{\bf 3}} \rput(56,82){{\bf 1}} \rput(8,82){{\bf 1}}
\rput(104,82){{\bf 4}}

\endpspicture
\caption{A semistandard $3$-ribbon tableaux with shape $(7,6,4,3,1)$, weight $(2,1,3,1)$ and spin 7.}
\label{fig:ribbontableaux}
\end{figure}

When $q=1$ the ribbon functions become products of usual Schur functions.
However, when the parameter $q$ is introduced, it is no longer obvious
that the functions $\g_\ll(X;q)$ are symmetric.  The main aim of this paper
will be to develop the theory of ribbon functions in the same way Schur
functions are studied in the ring of symmetric functions.  In particular,
we give:
\begin{itemize}
\item
A ribbon Pieri formula (Theorem \ref{thm:pieri}):
\[
h_k[\brac{1+q^{2}+\cdots+q^{2(n-1)}}X]\g_\nu(X;q) = \sum_\mu
q^{s(\mu/\nu)} \g_\mu(X;q).
\]
where the sum is over all $\mu$ such that $\mu/\nu$ is a
horizontal ribbon strip of size $k$.  The notation
$h_k[\brac{1+q^{2}+\cdots+q^{2(n-1)}}X]$ denotes a plethysm.
\item
A ribbon Murnagham-Nakayam-rule (Theorem \ref{thm:mnrule}):
\[
\brac{1+q^{2k}+\cdots+q^{2k(n-1)}}p_k\g_\nu(X;q) = \sum_\mu
\xx^{\mu/\nu}_k(q) \g_\mu(X;q).
\]
where $\xx^{\mu/\nu}_k(q)$ is combinatorially defined as an
alternating sum of spins over certain `border $n$-ribbon strips'
of size $k$.
\item
A Cauchy (and dual Cauchy) identity (Theorem \ref{thm:cauchy}):
\[
\sum_\ll \g_\ll(X;q)\g_\ll(Y;q) = \prod_{i,j} \prod_{k=0}^{n-1} \frac{1}{1-x_iy_jq^{2k}}
\]
where the sum is over all partitions $\ll$ with a fixed $n$-core. This formula had
been shown earlier \cite{Lam} for the case $n=2$.
\item
A $\C$-algebra isomorphism $\w_n: \Lambda(q) \rightarrow
\Lambda(q)$ (Theorem \ref{thm:inv}) satisfying
\[
\w_n(\g_\ll(X;q)) = \g_{\ll'}(X;q).
\]
Even the existence of a linear map with such a property is not obvious as the functions
$\g_\ll$ are not linearly independent.
\end{itemize}

It is well known that the corresponding formulae are important for Schur functions in
representation theory and algebraic geometry.

Much of the interest in the ribbon functions has been focused on the
$q$-Littlewood Richardson coefficients $c^\mu_{\ll}(q)$
of the expansion of $\g_\ll(X;q)$ in the Schur basis:
\[
\g_\ll(X;q) = \sum_\mu c^\mu_\ll(q) s_\ll(X).
\]
These are $q$-analogues of Littlewood Richardson coefficients.
Using results of Varagnolo and Vasserot \cite{VV}, Leclerc and
Thibon \cite{LT} have shown that these coefficients are parabolic
Kazhdan-Lusztig polynomials of type $A$.  Results of Kashiwara and
Tanisaki \cite{KT} then imply that they are polynomials in $q$
with non-negative coefficients.  Much interest has also developed
in connecting ribbon tableaux and the $q$-Littlewood Richardson
coefficients to rigged configurations and the generalised Kostka
polynomials defined by Kirillov and Shimozono \cite{KS}, Shimozono
and Weyman \cite{SW3}, Schilling and Warnaar \cite{SchW} and Shimozono \cite{Shi}. We
will focus mainly on the functions $\g_\ll(X;q)$ though clearly
our results imply relations between the $c^\mu_\ll(q)$. One
novelty is that we will study these functions even when the
$n$-core is non empty and more generally the $\g_{\ll/\mu}(X;q)$
for skew shapes $\ll/\mu$.

\medskip

To prove that the functions $\g_\ll(X;q)$ were symmetric Lascoux, Leclerc and Thibon
connected them to Fock space representation $\f$ of the quantum affine algebra $\uqsln$.
The crucial property of $\f$ is an action of a Heisenberg algebra $H$, commuting with the
action of $\uqsln$, discovered by Kashiwara, Miwa and Stern \cite{KMS}.  In particular,
they showed that as a $\uqsln \times H$-module, $\f$ decomposes as
\[
\f \cong V_{\Lambda_0} \otimes \CH
\]
where $V_{\Lambda_0}$ is the highest weight representation of $\uqsln$ with highest weight
$\Lambda_0$ and $\CH$ is the usual Fock space representation of the Heisenberg algebra.  We
will give a description of the action of the genrators $B_k$ of the $H$ on $\f$ in terms
of `border ribbon strips' by giving a combinatorial proof that the ribbon Murnagham-Nakayama
rule and the ribbon Pieri rules are formally equivalent (Theorem \ref{thm:mnviapieri}).

The connection between ribbon functions and the action of the
Heisenberg algebra is made explicit by showing (Theorem \ref{thm:map}) that the map $\m:\f \rightarrow \CH$ defined by
\[
\br{\ll} \mapsto \g_\ll
\]
is a map of $H$-modules, after identifying $\CH$ with the ring of symmetric functions $\Lambda(q)$
in the usual way.  The map $\m$ has the further remarkable property that it changes
certain linear maps into algebra maps, as follows.

Lascoux, Leclerc and Thibon \cite{LLT1} have constructed a global
basis of $\f$ which extends Kashiwara's global crystal basis of
$V_{\Lambda_0}$.  They defined a bar involution on $\f$ which
extends Kashiwara's involution on $V_{\Lambda_0}$.  Another
semi-linear involution, denoted $v \mapsto v'$ was also introduced
and further studied in \cite{LT} which satisfied the property
\[
\ip{\ov{u},v}=\ip{u',\ov{v'}}
\]
for $u,v \in \f$ and $\ip{\br{\ll},\br{\mu}} = \d_{\ll\mu}$ the
standard inner product on $\f$. We shall see that both involutions
become algebra isomorphisms under the map $\m$.  In particular the
`image' of the involution $v \mapsto v'$ is simply $\w_n$ (though
it is not immediately clear that such an image can be
well-defined) . We will also describe the image of certain vectors
in the global basis under the map $\m$.

Finally, we shall connect our study of ribbon functions to more
combinatorial aspects of ribbon tableaux. For the easiest non
trivial case of dominoes ($n=2$), we shall prove all our main
results using the domino insertion algorithms first studied by
Barbasch, Vogan and Garfinkle \cite{BV,Gar}.  Shimozono and White
\cite{SW} subsequently generalised domino insertion to the
semistandard case and also observed that the algorithm was
compatible with the spin statistic on domino tableaux.  In
\cite{Lam}, the author described a dual-domino Schensted algorithm
and observed both the Cauchy and dual-Cauchy identities.  Here, we
will prove the Pieri and Murnagham-Nakayama formulae using domino
insertion.

Shimozono and White \cite{SW2} have defined a ribbon-Schensted algorithm for $n > 2$ which is also compatible
with spin on ribbon tableaux.  As we shall discuss, this algorithm gives a combinatorial proof of the first ribbon
Pieri formula for $k=1$, but appears to be insufficient to prove either the Cauchy identity or the higher Pieri
rules.

The combinatorial approach to ribbon tableaux has been relegated to a secondary role in our presentation.  However,
it should be noted that the investigation of the Heisenberg algebra was inspired by empirical calculations with
domino and ribbon tableaux made while writing \cite{Lam}.  Throughout the paper we will use classical
symmetric function notation, however, most of our results could easily have been phrased in terms of the
Heisenberg algebra $H$.

\textbf{Organisation.} The article begins with two introductory
sections which give the notation we will use for tableaux and
symmetric functions. The main body of the paper is split into
three parts, which can be roughly described as being
representation theory, symmetric function theory and
combinatorics.  The reader interested mostly in the representation
theory will find that the first two parts can be read with nearly
no references to the last part. In Part \ref{part:rep}, we begin
by describing the action of $\uqsln$ on the Fock space $\f$.  The
details of this action will rarely be used in the paper, but we
present them for completeness. In Section \ref{sec:hei} we will
define the Heisenberg algebra and describe its action on both its
usual Fock space representation $\CH$ and on $\f$, as studied in
\cite{KMS, LLT, LT}. In Section \ref{sec:can} we will define the
global basis of $\f$ and the two involutions introduced by
Lascoux, Leclerc and Thibon.  While this section is important for
the overall understanding of the subject, it is logically
independent of most of the proofs in Part \ref{part:sym} which
mostly rely on the action of the Heisenberg algebra.  Only one new
result is present in Part \ref{part:rep}, a description of the
action of the generators $B_k$ of the Heisenberg algebra on $\f$
in terms of border ribbon strips. In Part \ref{part:sym}, we begin
by describing the initial properties of ribbon functions.  In
Section \ref{sec:canrib} we will relate the global basis of $\f$
to the ribbon functions.  In Section \ref{sec:mnrule}, we will
prove the ribbon Murnagham-Nakayama rule.  In Section
\ref{sec:map} we define and study the map $\m: \f \rightarrow
\Lambda(q)$.  In Section \ref{sec:pieri}, the Pieri rule is shown
modulo Theorem \ref{thm:mnviapieri} of Section
\ref{sec:mnviapieri}.  In Section \ref{sec:inv}, we introduce the
ribbon involution $\w_n$ and prove its main properties.  In
Section \ref{sec:cauchy}, we prove the Cauchy and dual-Cauchy
identities.  In Section \ref{sec:ip}, we describe a ribbon inner
product and study its relationship with another involution on
$\Lambda(q)$. In Section \ref{sec:skew}, we prove a `skew Cauchy
identity', and also define super ribbon tableaux and functions.
In Section \ref{sec:open}, we discuss some open problems. Part
\ref{part:comb} contains two rather separate sections.  Section
\ref{sec:ribinsertion} discusses the relationship between the
ribbon function formulae of Part \ref{part:sym} and ribbon
insertion algorithms.  Section \ref{sec:mnviapieri} contains a
standalone and purely combinatorial proof that the
Murnagham-Nakyama and Pieri rules are formally equivalent.

\textbf{Acknowledgements.}
This work is part of my dissertation written under the guidance of Richard Stanley.
I am indebted to him for suggesting the study of ribbon tableaux and
for providing me with assistance throughout.  I would also like to thank
Mark Shimozono and Ole Warnaar for pointing out a number of references.

\section*{Partitions and Tableaux} In this
section we give the notation and definitions we use for partitions
and ribbon tableaux.  A distinguished integer $n \geq 1$ will be
fixed throughout the whole.  When $n=1$, the reader may check that
we recover the classical theory of Schur functions.

A partition $\ll = (\ll_1 \geq \ll_2 \geq \cdots \geq \ll_l > 0)$
is a list of non-increasing integers.  We will call $l$ the length
of $\ll$, and denote it by $l(\ll)$.  We will say that $\ll$ is a
partition of $\ll_1+ \ll_2 + \ldots + \ll_l = |\ll|$ and write
$\ll \vdash |\ll|$.  We will typically use $\ll$, $\mu$, $\nu$ and
$\rho$ to denote partitions.  A composition $\a = (\a_1, \a_2,
\ldots, \a_l)$ is an ordered list of non-negative integers.  As
above, we will say that $\a$ is a composition of $|\a|=\a_1+\a_2+
\cdots +\a_l$.  Let $\ll$ and $\mu$ be partitions.  We will
generally not distinguish between a partition $\ll$ and its
corresponding Young diagram $D(\ll)$.  We will thus write $\ll
\subset \mu$ if $D(\ll) \subset D(\mu)$.  The skew shape $\ll/\mu$
is the set difference between the corresponding diagrams of $\ll$
and $\mu$ when $\mu \subset \ll$. The conjugate of a partition
$\ll$ obtained by changing rows to columns, is denoted $\ll'$.  We
will use the notation $m_k(\ll)$ to denote the number of parts of
$\ll$ equal to $k$.

\medskip
A skew shape $\ll/\mu$ is a horizontal strip if it contains at
most one square in each column.

A skew shape $\lambda/\mu$ is a border strip if it
 is connected, and does not contain any $2 \times 2$
square.  The height $h(b)$ of a border strip $b$ is the number of
rows in it, minus 1.  A border strip tableaux is a chain of
partitions
\[
\mu^0 \subset \mu^1 \subset \cdots \subset \mu^r
\]
such that each $\mu^{i+1}/\mu^{i}$ is a border strip.  The height
of a border strip tableaux $T$ is the sum of the heights of
its border strips.

When a border strip has $n$ squares for the distinguished (fixed)
integer $n$, we will call it a ribbon. The height of the ribbon
$r$ will then be called its spin $s(r)$. The reader should be
cautioned that in the literature the spin is usually defined as
half of this.

A semistandard tableaux of shape $\ll/\mu$ is a filling of each
square $(i,j)$ of the diagram $D(\ll/\mu)$ with a positive integer
such that the rows are non-decreasing and the columns are
increasing. The weight $w(T)$ of such a tableaux $T$ is the
composition $\a$ such that $\a_i$ is the number of occurrences of
$i$ in $T$.  The tableaux is standard if the numbers which occur
are exactly those of $[m]$ for some integer $m$.

\medskip

Let $\lambda$ be a partition.  Its $n$-core, obtained from
$\lambda$ by removal of $n$-ribbons (until we are no longer
able to), is denoted $\tilde{\lambda}$.
The $n$-quotient (see \cite{Mac}) of $\ll$ will be denoted
$(\ll^{(0)}, \ldots , \ll^{(n-1)})$.  We shall write $\p$ for the
set of partitions. We
will use $\p_\d$ to denote the set of partitions $\ll$ such that
$\tll = \d$ for an $n$-core $\d = \tilde{\d}$.

\begin{figure}[ht]
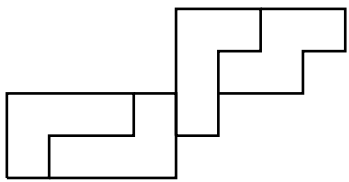

\pspicture(0,00)(130,80)

\psline(0,0)(0,32)(48,32)(48,16)(16,16)(16,0)(0,0)
\psline(16,0)(64,0)(64,32)(48,32)
\psline(64,16)(80,16)(80,48)(96,48)(96,64)(64,64)(64,16)
\psline(80,32)(112,32)(112,48)(128,48)(128,64)(96,64)

\endpspicture
\caption{A horizontal 4-ribbon strip with spin 5.}
\label{fig:horstrip}
\end{figure}

A ribbon tableaux $T$ of shape $\ll/\mu$ is a tiling of $\ll/\mu$
by $n$-ribbons and a filling of each ribbon with a positive
integer (see Figure \ref{fig:ribbontableaux}). If these numbers
are exactly those of $[m]$, for some $m$, then the tableaux is
called standard.  We will use the convention that a ribbon
tableaux of shape $\ll$ where $\tll \neq \emptyset$ is simply a
ribbon tableaux of shape $\ll/\tll$. A ribbon tableaux is
semistandard if for each $i$ \begin{enumerate}
\item
removing all ribbons labelled $j$ for $j > i$ gives a valid skew
shape $\ll_{\leq i}/\mu$ and,
\item
the subtableaux containing only the ribbons labelled $i$ form a
\emph{horizontal $n$-ribbon strip}.
\end{enumerate}
A horizontal $n$-ribbon strip is a skew shape tiled by ribbons
such that the topright-most square of every ribbon touches the
northern edge of the shape (see Figure \ref{fig:horstrip}). If
such a tiling exists, it is necessarily unique.

We will often think of a ribbon tableaux as a chain of partitions
\[
\tll = \mu^0 \subset \mu^1 \subset \cdots \subset \mu^r = \ll
\]
where each $\mu^{i+1}/\mu^{i}$ is a horizontal ribbon strip.  The
partitions $\mu^i$ here are not to be confused with the
$n$-quotient of $\mu$.

The spin $s(T)$ of a ribbon tableaux $T$ is the sum of the spins
of its ribbons.  The cospin $\text{cosp}(T)$ of a ribbon tableaux
$T$ is defined as $\text{cosp}(T) = \text{mspin}(sh(T)) - s(T)$,
where $\text{mspin}(sh(T))$ is the maximum spin of a ribbon
tableaux of the same shape as $T$. The weight $w(T)$ of a tableaux
is the composition counting the occurences of each value in $T$.

\medskip

Littlewood's $n$-quotient map (\cite{Lit}, see also \cite{SW1})
gives a weight preserving bijection between semistandard ribbon
tableaux $T$ of shape $\ll$ and $n$-tuples of semistandard Young
tableaux $\set{T^{(0)}, \ldots, T^{(n-1)}}$ of shapes $\ll^{(i)}$
respectively. Abusing language, we shall also refer to
$\set{T^{(0)}, \ldots, T^{(n-1)}}$ as the $n$-quotient of $T$.
Schilling, Shimozono and White \cite{SSW} have described the
cospin of a ribbon tableaux in terms of an inversion number of the
$n$-quotient. None of our proofs will require the use of the
$n$-quotient but occasionally we will comment on the $q=1$ case
for which the $n$-quotient will be important.

The $n$-quotient map can be described as follows.  A diagonal
$diag_d$ of a shape $\ll$ consists of all squares $(i,j)$ such
that $i-j = d$.  If we draw all diagonals of the form $diag_{dn}$
then each ribbon will intersect each such diagonal exactly once. A
ribbon's squares are linearly ordered from top right to bottom
left. Suppose the diagonal $diag_{dn}$ intersects a ribbon $r$ at
the $k^{th}$ square from the top right.  Then the ribbon $r$ is
sent under the $n$-quotient map to a square in the diagonal
$diag_d$ of $\ll^{(k)}$. The numbers in the ribbon tableaux of
Figure \ref{fig:ribbontableaux} have been placed along the
diagonals $diag_{dn}$.  Figure \ref{fig:quotient} shows its
$3$-quotient.

\begin{figure}[ht]
\pspicture(-20,40)(200,80) \rput(-20,55){$T^{(0)} = $}
\abox(0,4){2} \abox(1,4){3} \abox(2,4){4}
\rput(102,55){$T^{(1)} = $} \abox(8,4){1} \abox(9,4){1}
\abox(8,3){3} \abox(9,3){3} \rput(220,55){$T^{(2)} =
\emptyset $}
\endpspicture
\caption{The $3$-quotient of the ribbon tableaux $T$ of Figure
\ref{fig:ribbontableaux}.} \label{fig:quotient}
\end{figure}

A horizontal ribbon strip can be described in terms of the
$n$-quotient as a union of horizontal strips for each tableaux of
the $n$-quotient.

\section*{Symmetric Functions}
In this section we briefly review some standard notation in
symmetric function theory.  The reader is referred to \cite{Mac}
for further details.

Let $\Lambda_\zz$ denote the ring of symmetric functions with
coefficients in $\zz$.  Recall that $\Lambda_\zz$ has a
distinguished integral basis $s_\lambda$ known as the Schur
functions.  Nearly all the results of this paper can be stated in
$\Lambda_\zz[q]$, but some intermediate steps may require working
in $\Lambda = \Lambda_\C$ so we will use that as our symmetric
function ring from now on.  We will write $\Lambda(q)$ for
$\Lambda_\C \otimes_\C \C(q)$.

It is well known that the Schur functions $s_\ll$ are orthogonal
with respect to a natural inner product $\langle \, , \, \rangle$
on $\Lambda$ and are unique up to signed permutation.  We will
denote the homogeneous, elementary, monomial and power sum
symmetric functions by $h_\ll$, $e_\ll$, $m_\ll$ and $p_\ll$
respectively. Recall that we have $\ip{h_\ll,m_\mu} = \d_{\ll\mu}$
and $\ip{p_\ll,p_\mu} = z_\ll \d_{\ll\mu}$ where $z_\ll =
1^{m_1(\ll)} m_1(\ll)! 2^{m_2(\ll)} m_2(\ll)! \cdots$.  Each of
$\set{p_i}$, $\set{e_i}$ and $\set{h_i}$ generate $\Lambda$.  We will
write $X$ to mean $(x_1,x_2,\ldots)$.  Thus $s_\ll(X) = s_\ll(x_1,x_2,\ldots)$.

Recall that the Kostka matrix $K_{\ll\mu}$ is defined as
\[
s_\ll = \sum_\mu K_{\ll\mu} m_\mu.
\]
We will denote the inverse Kostka matrix by $\kappa_{\ll\mu}$:
\[
m_\mu = \sum_\ll \kappa_{\ll\mu} s_\ll.
\]

Let $f \in \Lambda$.  We will recall the definition of the
plethysm $g \mapsto g[f]$.  Write $g = \sum_\ll c_\ll p_\ll$. Then
we have
\[
g[f] = \sum_\ll c_\ll \prod_{i=1}^{l(\ll)} f(x_1^{\ll_i},
x_2^{\ll_i}, \ldots).
\]
Thus the plethysm by $f$ is the (unique) algebra isomorphism of
$\Lambda$ which sends $p_k \mapsto f(x_1^k, x_2^k, \ldots)$.  When
$f(x_1,x_2,\ldots; q) \in \Lambda(q)$ for a distinguished element
$q$, we define the plethysm as $p_k \mapsto f(x_1^k, x_2^k,
\ldots;q^k)$.

For example, the plethysm by $(1+q)p_1$ is given by sending
\[
p_k \mapsto (1+q^k) p_k
\]
and extending to an algebra isomorphism $\Lambda(q) \rightarrow
\Lambda(q)$.  In such situations we will write $f[(1+q)X]$ for
$f[(1+q)p_1]$.

We will be particularly concerned with the plethysm given by
$(1+q^2+ \cdots +q^{2n-2})p_1$.  We will use $\plqn$ to denote the
map $\Lambda(q) \rightarrow \Lambda(q)$ given by $f \mapsto
f[(1+q^2+\cdots+q^{2n-2})X]$.

\part{The Fock Space of $\uqsln$ and the Heisenberg Algebra}
\label{part:rep}
\section{The Fock Space representation $\f$ of $\uqsln$}
\label{sec:uqsln}
In this section we introduce the quantum affine algebra $\uqsln$
and its ($q$-deformed) Fock Space representation $\f$. We will
only be using the action of $\uqsln$ to define the canonical
basis, but we include the details for completeness.  A concise
introduction to the material of this section can be found in
\cite{Lec}.  Throughout $q$ can be thought of as either a formal
parameter of as a generic complex number (not equal to a root of
unity).

We denote by ${\mathfrak h}$ the `Cartan' subalgebra of
$\hat{{\mathfrak sl_n}}$ spanned over $\C$ by the basis
$\set{h_0,h_1,\ldots,h_{n-1},D}$.  The dual basis will be spanned
by $\set{\Lambda_0,\Lambda_1,\ldots,\Lambda_{n-1},\d}.$  We set
$\a_i = 2\Lambda_i - \Lambda_{i-1} -\Lambda_{i+1}$ for $i \in
\set{1,2,\ldots,n-2}$, and $\a_0 = \Lambda_0 - \Lambda_{n-1} -
\Lambda_1 + \d$ and $\a_{n-1} = - \Lambda_{n-2} + \Lambda_{n-1} -
\Lambda_0$. The generalised Cartan matrix $[\ip{\a_i,h_j}]$ will
be denoted $a_{ij}$.  Set $P^{\vee}= (\oplus_{i=0}^{n-1}\zz h_i)
\oplus \zz D$.

The algebra $\uqsln$ is the associative algebra over $\C(q)$
generated by elements $e_i$, $f_i$ for $0 \leq i \leq n-1$, and
$q^{h}$ for $h \in P^{\vee}$ satisfying the following relations:
\[
q^{h}q^{h'} = q^{h+h'},
\]
\[
q^{h_i}e_jq^{-h_i} = q^{a_{ij}}e_j, \;\;\; q^{h_i}f_jq^{-h_i} =
q^{-a_{ij}}f_j,
\]
\[
q^De_iq^{-D} = \d_{i0}q^{-1}e_0, \;\;\; q^Df_iq^{-D} =
\d_{i0}qf_0,
\]
\[
[e_i,f_j] = \d_{ij} \frac{q^{h_i}-q^{-h_i}}{q-q^{-1}},
\]
\[
\sum_{k=0}^{1-a_{ij}} (-1)^k \left[ \begin{array}{c} 1-a_{ij} \\ k
\end{array} \right] e_i^{1-a_{ij}-k}e_je_i^k = 0 \;\,\,\, \mbox{$(i \neq
j)$}, \] \[ \sum_{k=0}^{1-a_{ij}} (-1)^k \left[ \begin{array}{c}
1-a_{ij} \\ k
\end{array} \right] f_i^{1-a_{ij}-k}f_jf_i^k = 0 \;\,\,\,\mbox{$(i \neq
j)$}. \]

We have used the standard notation
\[
[k] = \frac{q^k-q^{-k}}{q-q^{-1}}, \;\;\; [k]! =
[k][k-1]\cdots[1],
\]
and
\[
\left[\begin{array}{c}n \\k \end{array}\right] =
\frac{[n]!}{[n-k]![k]!}.
\]

\medskip
The Fock Space $\f$ is an infinite dimensional vector space over
$\C(q)$ spanned by a countable basis $\br{\ll}$ indexed by $\ll
\in \p$.  We follow the terminology of \cite{LLT}.

There is an action of the quantum affine algebra $\uqsln$ on $\f$
due to Hayashi \cite{Hay} which was formulated essentially as
follows by Misra and Miwa \cite{MM}.

Recall that a cell $(i,j)$ has content given by $c(i,j) = i-j$.
Its residue $p(i,j) \in \set{0,1,\ldots,n-1}$ is then $i-j$ mod
$n$. We call $(i,j)$ an indent $k$-node of $\ll$ if $p(i,j)=k$ and
$\ll \cup (i,j)$ is a valid Young diagram.  We make the analogous
definition for a removable $i$-node.

\medskip
Let $i \in \set{0,1,\ldots,n-1}$ and $\mu = \ll \cup \d$ for an
indent $i$-node $\d$ of $\ll$.  Now set
\par
  $N_i(\ll) = \#\{$ indent $i$-nodes of $\ll\} - \#\{$ removable $i$-nodes
 of $\ll\}$.
\par $N_i^l(\ll,\mu) = \#\{$ indent $i$-nodes of $\ll$ to the left of
$\d$ (not counting $\d$) $\} - \#\{$ removable $i$-nodes of $\ll$
to the left of $\d\}$.
\par
$N_i^r(\ll,\mu) = \#\{$ indent $i$-nodes of $\ll$ to the right of
$\d$ (not counting $\d$) $\} - \#\{$ removable $i$-nodes of $\ll$
to the left of $\d\}$.
\par
$N^0(\ll) = \#\{$ 0 nodes of $\ll\}$.

\medskip
Then we have the following theorem.
\begin{thm}
The following formulae define an action of the quantum affine
algebra $\uqsln$ on $\f$:
\begin{align*}
&q^{h_i}\br{\ll} = q^{N_i(\ll)}\br{\ll}, \;\mbox{for each $i \in \set{0,1,\ldots,n-1}$}, \\
&q^{D}\br{\ll} = q^{N^0(\ll)}\br{\ll}, \\
f_i\br{\ll} &= \sum_\mu q^{N_i^r(\ll,\mu)}\br{\mu}, \;
\mbox{summed
over all $\mu$ such that $\mu/\ll$ has residue $i$}, \\
e_i\br{\ll} &= \sum_\mu q^{-N_i^l(\ll,\mu)}\br{\mu}, \;
\mbox{summed
over all $\mu$ such that $\ll/\mu$ has residue $i$}. \\
\end{align*}
\end{thm}

The $\uqsln$-submodule of $\f$ generated by the vector $\br{0}$ is
easily seen to be the irreducible highest weight module with
highest weight $\Lambda_0$, which we will denote $V_{\Lambda_0}$.

\section{The Action of the Heisenberg Algebra}
\label{sec:hei}
This action of the Heisenberg Algebra on $\f$ will be essential
for our study of ribbon functions in Part \ref{part:sym}.

The Heisenberg Algebra $H$ will be the associative algebra with 1
generated over $\C(q)$ by a countable set of generators $B_k: k
\in \zz - \set{0}$ satisfying \be{eq:hei} [B_k,B_l] = la_l(q)
\d_{k,-l} \ee for some elements $a_l(q) \in \C(q)$ satisfying
$a_l(q) = a_{-l}(q)$.  (Often the element 1 is called the central
element and denoted $c$, but we will not need this generality).
\medskip The Fock Space representation $\CH$ of
$H$ is the polynomial algebra
\[
\CH \cong \C(q)[B_{-1},B_{-2}, \ldots].\]  The elements $B_{-k}$
for $k \geq 1$ act by multiplication on $\CH$. The action of $B_k$
for $k \geq 1$ is given by (\ref{eq:hei}) and the relation
\[
B_k\cdot 1 = 0 \; \mbox{for $k \geq 1$.}
\]

One common explicit construction of $\CH$ is given by
$\Lambda(q)$.  We may identify $B_k$ as the following operators:
\[
B_{-k}: f \mapsto a_k(q)p_k\cdot f \; \mbox{for $k \geq 1$}
\]
and
\[
B_k: f \mapsto k \frac{\partial}{\partial p_k} f \; \mbox{for $k
\geq 1$.}
\]
Under this identification, the operators $B_k$ have degree $-k$.
The reason for splitting $a_k(q)$ and $k$ is because $p_k$ and
$k\frac{\partial}{\partial p_k}$ are adjoint operators under the
usual inner product $\langle \, , \, \rangle$ on $\Lambda(q)$ (see
for example \cite{Mac}).

A standard lemma that we shall need later is
\begin{lem}
\label{lem:bkcommute} Let $k \geq 1$ be an integer and $\ll$ be a
partition. Then
\[
B_k B_{-\ll} = ka_k(q)m_k(\ll)B_{-\mu} + B B_k
\]
where $m_k(\ll)$ is the number of parts of $\ll$ equal to $k$ and
$B$ is some element in $H$ and $\mu$ is $\ll$ with one less part
equal to $k$. If $m_k(\ll) = 0 $ to begin with then the first term
is just 0.
\end{lem}
\begin{proof}
We may commute $B_k$ with $B_{-\ll_i}$ immediately for parts
$\ll_i \neq k$.  For each part equal to $k$, using the relation
$[B_{-k},B_{k}] = ka_k(q)$ introduces one term of the form
$ka_k(q)B_{-\mu}$.
\end{proof}

Kashiwara, Miwa and Stern \cite{KMS} have defined an action of the
the affine Hecke algebra $\hat{H}_N$ on the tensor product $V(z)^{\otimes N}$
of evaluation modules.
As $N \rightarrow \infty$, we obtain an action of the center $Z(\hat{H}_N)$ as a copy
of the Heisenberg algebra $H$ on $\f$, commuting with the action of $\uqsln$.
The operators $B_k$ are given as the infinite power
sums
\[
B_k = \sum_{i=1}^{\infty} Y_i^{-k}
\]
in terms of the certain generators $Y_i$ (in the $\set{T_i,Y_j}$
presentation) of $\hat{H}_N$.

The following theorems are due to Kashiwara, Miwa and Stern
\cite{KMS}.

\begin{thm}
\label{thm:heisenberg} The operators $B_k$ commute with the action
of the quantum affine algebra $\uqsln$.  They satisfy the
relations
\[
[B_k,B_l] = k \frac{1-q^{-2nk}}{1-q^{-2k}} \d_{k,-l}
\]
and generate a copy of the Heisenberg algebra.
\end{thm}
We shall see later in Section \ref{sec:ribinsertion} that the
factor $\brac{\frac{1-q^{-2nk}}{1-q^{-2k}}}$ can be given a
combinatorial explanation in terms of ribbon insertion.

\begin{thm}
The Fock space $\f$, regarded as a representation of $\uqsln
\otimes U(H)$ decomposes as the tensor product
\[
\f \simeq V_{\Lambda_0} \otimes \CH
\]
where $\CH$ is the Fock space of the Heisenberg algebra $H$ and
$V_{\Lambda_0}$ is the highest weight representation with highest
weight $\Lambda_0$.
\end{thm}

We write $B_\a = B_{\a_l} B_{\a_{l-1}} \cdots B_{\a_1}$ for a
composition $\a$. Similarly, $B_{-\a}$ denotes $B_{-\a_l} \cdots
B_{-\a_1}$.

Lascoux, Leclerc and Thibon \cite{LLT} have described the action
of certain elements $\u_k$, $\tu_k$, $\v_k$, and $\tv_k$ of $H$ on
the Fock Space $\f$, in terms of ribbon tableaux.

In terms of the elements $Y_i^{\pm 1}$ of the affine Hecke
algebra, we have
\begin{align*}
\u_k &= h_k(Y_1,Y_2,\ldots), \\
\v_k &= h_k(Y_1^{-1},Y_2^{-1},\ldots), \\
\tu_k &= e_k(Y_1,Y_2,\ldots), \\
\tv_k &= e_k(Y_1^{-1},Y_2^{-1},\ldots), \\
\end{align*}
where the $e_k$ and $h_k$ are the elementary and homogeneous
symmetric functions. To avoid mention of the elements $Y_i$ one
may write them as
\[
\u_k = \sum_\ll b_{k,\ll}B_\ll
\]
where the coefficients $b_{k,\ll}$ are given by the expansion
\[
h_k = \sum b_{k,\ll} p_\ll
\]
in the ring $\Lambda$ of symmetric functions and $p_\ll$ are the
power sum symmetric functions.

\begin{prop}
\label{prop:hor} The elements $\u_k$, $\tu_k$, $\v_k$, and $\tv_k$
of $\h$ act on $\f$ as linear operators defined by
\[
\v_k \br{\ll} = \sum_\mu (-q)^{-s(\mu/\ll)} \br{\mu},
\]
where the sum is over all $\mu$ such that $\mu/\ll$ is a
horizontal $n$-ribbon strip of size $k$.  Similarly,
\[
\u_k \br{\ll} = \sum_\nu (-q)^{-s(\ll/\nu)} \br{\nu},
\]
summed over all $\nu$ such that $\ll/\nu$ is a horizontal
$n$-ribbon strip of size $k$.  The formulae for $\tu_k$ and
$\tv_k$ are exactly analagous, with horizontal ribbon strips
replaced by vertical ribbon strips.
\end{prop}

We will write $\v_\a$ for $\v_{\a_l}\cdots\v_{\a_1}$ and similarly
for $\u_\a$, $\tv_\a$ and $\tu_\a$.  Thus
\[
\v_\a \br{\mu} = \sum_T q^{s(T)}\br{\ll}
\]
summed over all semistandard ribbon tableaux of shape $\ll/\mu$
and weight $\a$.

In Part $\ref{part:comb}$ we will show (Theorem \ref{thm:mnviapieri})
that the Pieri style rule
of Proposition \ref{prop:hor} formally implies a Murnagham-Nakayama
style rule for the the operators $B_k$ and $B_{-k}$.
The precise statement and proof
of Theorem \ref{thm:mnviapieri} has been deferred until the end as
its proof is completely combinatorial and unrelated to the Fock
space.

\begin{prop}
\label{prop:bkaction} The linear operators $B_{-k}$ for $k \geq 1$
act on $\f$ as
\[
B_{-k} \br{\ll} = \sum_\mu \xx_{k}^{\mu/\ll}(-q^{-1}) \br{\mu}
\]
where $\xx_{k}^{\mu/\ll}(q)$ is given by
\[
\xx_{k}^{\mu/\ll}(q) = \sum_S (-1)^{h(S)}q^{s(S)}
\]
summed over all border ribbon strip tilings $S$ of $\mu/\ll$.
Similarly,
\[
B_{k} \br{\ll} = \sum_\mu \xx_{k}^{\ll/\mu}(-q^{-1}) \br{\mu}.
\]
\end{prop}

\begin{proof}
This is an immediate consequence of Proposition \ref{prop:hor} and
Theorem \ref{thm:mnviapieri}.
\end{proof}
Border ribbon strips will be defined combinatorially later (Definition
\ref{def:borderribbon}) and are ribbon analogues of usual border
strips.

It follows immediately from the above propositions and Theorem
\ref{thm:heisenberg} that the sets $\set{\v_\ll \br{0}}_{\ll \in
\p}$, $\set{\u_\ll \br{0}}_{\ll \in \p}$ and $\set{B_\ll
\br{0}}_{\ll \in \p}$ form bases of the space of highest weight
vectors of $\uqsln$ in $\f$.

\section{Global Bases of $\f$}
\label{sec:can} We first define a bar involution $v \mapsto
\ov{v}$ on $\f$ following Leclerc and Thibon \cite{LT,LT1}.  This
involution restricted to $V_{\Lambda_0}$ (the $\uqsln$ submodule
with highest weight vector $\br{0}$) agrees with
Kashiwara's involution \cite{Kas}.

\begin{prop}
\label{prop:bar} There exists a unique semi-linear map ${}^{-}: \f
\rightarrow \f$ satisfying
\begin{align*}
\ov{qv} &= q^{-1}\ov{v}, \\
\ov{f_i\cdot v} &= f_i \cdot \ov{v}, \\
\ov{e_i\cdot v} &= e_i \cdot \ov{v}, \\
\ov{B_{-k} \cdot v} &= B_{-k} \cdot \ov{v}, \\
\ov{B_k \cdot v} &= q^{2(n-1)k}B_k \cdot \ov{v}.\\
\end{align*}
\end{prop}

The Fock space $\f$ has a natural order `$<$' defined by $\br{\ll}
< \br{\mu}$ if and only if $\ll \prec \mu$ in dominance order.
Leclerc and Thibon show that ${}^{-}$ is triangular with respect
to the basis $\br{\ll}$ and conclude that the global basis of the following
theorem exists.

\begin{thm}
\label{thm:canexist} There exist unique vectors $G_{\ll} \in \f$
for $\ll \in \p$ satisfying:
\[
\ov{G_{\ll}} = G_{\ll}
\]
and
\[
G_{\ll}\equiv \br{\ll} \;\mbox{mod $q^{-1}{\mathcal L}^-$}
\]
where ${\mathcal L}^-$ is the $\zz[q^{-1}]$ submodule of $\f$
spanned by $\br{\ll}$.
\end{thm}
When we restrict this to the $\uqsln$ submodule $V_{\Lambda_0}$ of $\f$
generated by $\br{0}$, the $G_{\ll}$ is essentially the global
upper crystal basis of $V_{\Lambda_0}$ (see \cite{LLT1,Kas1}). This follows from
the fact that the bar involution agrees with Kashiwara's
involution when restricted to $V_{\Lambda_0}$.

Some of the $G_{\ll}$ are especially easy to describe in terms of
the action of the Heisenberg algebra.  In particular, in analogy
with Steinberg's tensor product theorem, Leclerc and Thibon
\cite{LT} show that $G_{\mu} = S_\ll\br{\nu}$ for an $n$-regular
partition $\nu$ and $\mu = n\ll + \nu$. We will only need the
following special case.

\begin{prop}
\label{prop:can} Let $S_\ll = \sum_\mu \Chi^\ll_\mu B_{-\ll}$.
Thus in terms of the elements $Y_i$ we have $S_\ll =
s_\ll(Y_1^{-1}, Y_2^{-1}, \ldots)$.  Then we have
\[
G_{n\ll} = S_\ll \br{0}.
\]
\end{prop}
\begin{proof}
It is clear that $\ov{\br{0}} = \br{0}$.  Thus by Proposition
\ref{prop:bar} we have $\ov{S_\ll\br{0}} = S_\ll\br{0}$.  By the
definition of $G_{n\ll}$ it suffices to show that $S_\ll\br{0}
\equiv \br{n\ll}$ mod $q^{-1}{\mathcal L}^-$.

By Proposition \ref{prop:hor} we know that
\[
\v_\a \br{0} \equiv \sum_T \br{sh(T)} \;\; \mbox{mod
$q^{-1}{\mathcal L}^-$}
\]
where the sum is over all ribbon tableaux of spin 0 and weight
$\a$.  It is clear that these are in bijection with usual
semistandard Young tableaux of weight $\a$.  Thus
\[
\v_\a \br{0} \equiv \sum_\ll K_{\ll \a} \br{n\ll} \;\; \mbox{mod
$q^{-1}{\mathcal L}^-$}.
\]
Comparing with $\v_\a = \sum_\ll K_{\ll\a} S_\ll$ we see that
\[
S_\ll \br{0} \equiv \br{n\ll} \;\; \mbox{mod $q^{-1}{\mathcal
L}^-$}
\]
which completes the proof.
\end{proof}
It follows immediately that $\set{G_{n\ll}}$ form a basis of the
space of highest weight vectors of the action of $\uqsln$ on $\f$.

\begin{remark}
\label{rem:can}
\smallskip

\begin{enumerate}
\item
Let
\[
G_{\ll} = \sum_\mu l_{\ll,\mu}(-q^{-1}) \br{\mu}.
\]
Varagnolo and Vasserot \cite{VV} have shown that $l_{\ll,\mu}$ is
a parabolic Kazhdan-Lusztig polynomial for the affine Hecke
algebra of type $A$.  These polynomials were introduced by Deodhar
\cite{D1,D2} and shown to have non-negative coefficients by
Kashiwara and Tanisaki \cite{KT}.
\item
We have not mentioned the lower global basis $G^+$ of $\f$ as they
are less related to the ribbon functions we will be studying. The
$G^+_\ll$ are defined in a similar way to the $G_\ll$ by
\begin{align*}
G_{\ll}^+ &= \ov{G_\ll^+}\\
G_{\ll}^+ &\equiv \br{\ll} \;\mbox{mod $q{\mathcal L}^+$}
\end{align*}
where ${\mathcal L}^+$ is the $\zz[q]$-submodule spanned by
$\br{\ll}$.
\end{enumerate}

\end{remark}

Finally, we will be needing a semi-linear involution $v \mapsto
v'$ on $\f$.  This is defined by $q' = q^{-1}$ and
\[
\br{\ll} \mapsto \br{\ll'}.
\]

Then we have \cite[Proposition 7.10]{LT}
\begin{prop}
\label{prop:semilinear} For all $v \in \f$ and compositions $\a$
satisfying $|\a| = k$ we have
\begin{align*}
&(e_iu)' = q^{h_{-i}-1}e_{-i}u', &(f_iu)' = q^{-h_{-i}-1}f_{-i}u',
\\ &(\v_\beta u)' = (-q)^{(n-1)k} \tv_\beta u', &(\u_\beta u)' = (-q)^{(n-1)k} \tu_\beta u'.  \\
\end{align*}
\end{prop}

One immediate consequence of this and Proposition \ref{prop:can}
is that
\[
\brac{G_{n\ll}}' = (-q)^{(n-1)k}G_{n\ll'}.
\]

\medskip
We shall see later an explicit connection between the global basis
and ribbon tableaux.  This should come as no surprise: the $\v_\a$
are described combinatorially in terms of ribbon tableaux and
$S_\ll$ can be expressed in terms of the $\v_\a$ (via the inverse
Kostka matrix).

\part{Ribbon Functions}
\label{part:sym}

\section{Definitions and Initial Properties}
\label{sec:not}
We will now define the central objects of this paper as introduced
by Lascoux, Leclerc and Thibon in \cite{LLT}.

\begin{definition}
Let $\ll/\mu$ be a skew partition, tileable by $n$-ribbons. Define
the symmetric functions $\g_{\ll/\mu} \in \Lambda(q)$ as:
\[
\g_{\ll/\mu}(X;q) = \sum_T q^{s(T)}{\mathbf x}^{w(T)}
\]
where the sum is over all semistandard ribbon tableaux $T$ of
shape $\ll/\mu$ and ${\mathbf x}^\a = x_1^{\a_1} x_2^{\a_2}
\cdots$. When $\ll$ is a partition with non-empty $n$-core, we
write $\g_\ll$ for $\g_{\ll/\tll}$.  These functions will be
loosely called \emph{ribbon functions}.
\end{definition}

The fact that the functions $\g_{\ll/\mu}$ are symmetric is not
obvious from the combinatorial definition.  However, using the
action of the Heisenberg algebra on the Fock space $\f$, the proof
is immediate (\cite{LLT}) and reproduced below.

\begin{definition}
Let $\ll/\mu$ be a skew shape tileable by $n$-ribbons.  Then
define
\[
\kk_{\ll/\mu,\a}(q) = \sum_T q^{s(T)},
\]
the spin generating function of all semistandard ribbon tableaux
$T$ of shape $\ll/\mu$ and weight $\a$.
Similarly let
\[
\l_{\ll/\mu,\a}(q) = \sum_T q^{s(T)}
\]
summed over all column semistandard ribbon tableaux of shape $\ll/\mu$
and weight $\a$.  A ribbon tableaux is column semistandard
if its conjugate is semistandard.
\end{definition}

Thus
\[
\g_{\ll/\mu}(X;q) = \sum_\a \kk_{\ll/\mu,\a}(q) {\mathbf x}^\a.
\]

\begin{thm}
\label{thm:sym}The functions $\g_{\ll/\mu}(X;q)$ are symmetric
functions.
\end{thm}
\begin{proof}
A semistandard ribbon tableaux can be expressed as a chain of
partitions differing by horizontal ribbon strips.  Thus
\[
\v_\a \br{\mu} = \sum_\nu \kk_{\nu/\mu, \a}(-q^{-1}) \br{\nu}.
\]
But if $\beta$ is a permutation of $\a$, then $\v_\a = \v_\b$
since the $\v_k$ commute.  This shows that
\[
\kk_{\ll/\mu,\a}(q) = \kk_{\ll/\mu,\b}(q),
\]
after equating coefficients of $\ll$.
\end{proof}

We will also need the following definition of a \emph{border strip
ribbon tableaux}.
\begin{definition}
\label{def:borderribbon} A \emph{border ribbon strip} $T$ is a
connected skew shape $\ll/\mu$ with a distinguished tiling by
disjoint non-empty horizontal ribbon strips $T_1, \ldots, T_a$
such that the diagram $T_{+i} = \cup_{i \leq j} T_i$ is a valid
skew shape for every $i$ and for each connected component $C$ of
$T_i$ we have
\begin{enumerate}
\item
The shape of $C \cup T_{i-1}$ is not a horizontal ribbon strip.
Thus $C$ has to `touch' $T_{i-1}$ `from below'.
\item
No sub horizontal ribbon strip $C'$ of $C$ which can be added to $T_{i-1}$ satisfies the above
property.  Since $C$ is connected, this is equivalent to saying that only the rightmost ribbon of $C$ touches 
$T_{i-1}$.
\end{enumerate}
We further require that $T_1$ is connected.  The height $h(T_i)$
of the horizontal ribbon strip $T_i$ is the number of its
components. The height $h(T)$ of the border ribbon strip is
defined as $h(T) = \left(\sum_i h(T_i)\right) - 1$. The size of
the border ribbon strip $T$ is then the total number of ribbons in
$\cup_i T_i$.  A border ribbon strip tableaux is a chain $T =
\ll_0 \subset \ll_1 \cdots \subset \ll_r$ of shapes such that
$\ll_i/\ll_{i-1}$ has been given the structure of a border ribbon
strip.  The type of $T = \set{\ll_i}$ is then the composition $\a$
with $\a_i$ equal to the size of $\ll_i/\ll_{i-1}$.

Define $\xx^{\mu/\ll}_{\nu}$ as
\[
\xx^{\mu/\ll}_{\nu}(q) = \sum_T (-1)^{h(T)}q^{s(T)}
\]
summed over all border ribbon strip tableaux of shape $\mu/\ll$
and type $\nu$.

\end{definition}
Note that this definition reduces to the usual definition of a
border strip and border strip tableaux
when $n=1$, in which case all the horizontal strips $T_i$ are actually connected.

\begin{figure}[ht]
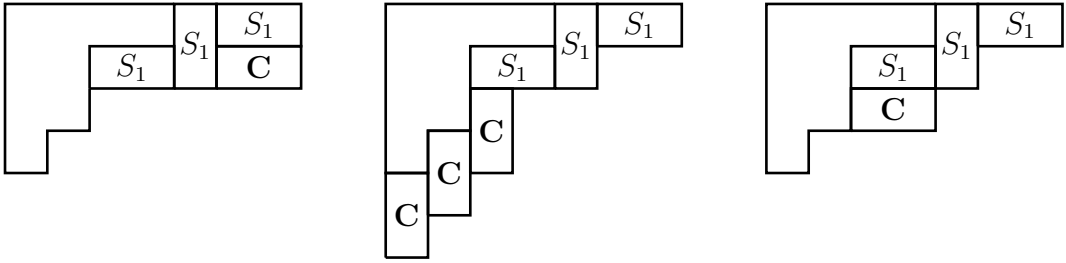

\pspicture(40,10)(450,120)

\hdom(5,5){$S_1$}
\vdom(7,5){$S_1$}
\hdom(8,6){$S_1$}
\hdom(8,5){C}
\psline(80,80)(80,64)(64,64)(64,48)(48,48)(48,112)(112,112)

\vdom(12,1){C}
\vdom(13,2){C}
\vdom(14,3){C}
\hdom(14,5){$S_1$}
\vdom(16,5){$S_1$}
\hdom(17,6){$S_1$}
\psline(224,80)(224,64)(208,64)(208,48)(192,48)(192,112)(256,112)

\hdom(23,5){$S_1$}
\vdom(25,5){$S_1$}
\hdom(26,6){$S_1$}
\hdom(23,4){C}
\psline(368,80)(368,64)(352,64)(352,48)(336,48)(336,112)(400,112)

\endpspicture
\caption{Connected horizontal strips $C$ which can be added  to $S_1 =(7,5,2,1)/(4,2,2,1)$ to form a border ribbon strip.  The resulting
border ribbon strips all have height 1.}
\label{fig:borderribbonstrip}
\end{figure}

\begin{example}
Let $n=2$ and $\ll = (4,2,2,1)$.  Suppose $S$ is a border ribbon strip such that $S_1$ has shape
$(7,5,2,1)/(4,2,2,1)$, and thus it has size 3 and spin 1.  We will now determine all the possible horizontal ribbon strips
which may form $S_2$.  It suffices to find the possible connected components that may be added.  The domino $(9,5,2,1)/(7,5,2,1)$
may not be added since its union with $S_1$ is a horizontal ribbon strip, violating the conditions of the definition.  The
domino strip $(8,8,2,1)/(7,5,2,1)$ is not allowed since the domino $(8,8,2,1)/(7,7,2,1)$ can be removed and we still obtain a 
strip which touches $S_1$. 

The legitimate connected horizontal ribbon strips $C$ which can be added are $(7,7,2,1)/(7,5,2,1)$, $(7,5,4,1)/(7,5,2,1)$ and $(7,5,3,3,2,1)/(7,5,2,1)$ as shown in 
Figure \ref{fig:borderribbonstrip}.  Thus assuming $S_2$ is non-empty, there are 5 choices for $S_2$, corresponding to taking some 
compatible combination of the three connected horizontal ribbon strips above.
\end{example}

\begin{example}
As before let $n=2$.  We will calculate $\xx^{\ll/\mu}_5(q)$ for $\ll = (5,5,2)$ and $\mu = (2)$.
The relevant border ribbon strips $S$ are (successive differences of the following chains denote the $S_i$)
\begin{itemize}
\item
$(2) \subset (5,5,2)$ with height 0 and spin 5, 
\item
$(2) \subset (5,3,2) \subset (5,5,2)$ with height 1 and spin 3, 
\item
$(2) \subset (5,5) \subset (5,5,2)$ with height 1 and spin 3,
\item
$(2) \subset (5,3) \subset (5,5,2)$ with height 2 and spin 1. 
\end{itemize}
Thus
\[
\xx^{\ll/\mu}_5(q) = q^5 - 2q^3 + q.
\]

\end{example}

When $q=1$, the ribbon functions become products of Schur
functions:
\[
\g_\ll(X;1) = s_{\ll^{(0)}}s_{\ll^{(1)}}\cdots s_{\ll^{(n-1)}}.
\]
This is a consequence of Littlewood's $n$-quotient map.  In fact,
up to sign, $\g_\ll(X;1)$ is essentially $\phi_n(s_\ll)$ where
$\phi_n$ is the adjoint operator to taking the plethysm by $p_n$.
More generally, $\g_{\ll/\mu}$ reduces to a product of skew Schur
functions at $q=1$:
\[
\g_{\ll/\mu}(X;1) = s_{\ll^{(0)}/\mu^{(0)}}s_{\ll^{(1)}/\mu^{(1)}}\cdots s_{\ll^{(n-1)}/\mu^{(n-1)}}.
\]

\begin{remark}
Another set of symmetric functions $\h_\ll(X;q)$ are defined
(\cite{LLT}) by
\[
\h_\ll(X;q) = \g_{n\ll}(X;q).
\]
It is not hard to see that
\[
\h_\ll(X;0) = s_\ll
\]
and that
\[
\h_\ll(X;1) = s_\ll + \sum_{\mu \prec \ll}s_\mu
\]
where $\prec$ denotes the usual dominance order on partitions.
Thus the functions $\h_\ll(X;q)$ form a basis of $\Lambda(q)$ over
$\C(q)$.  In \cite{LLT} it is shown that the cospin versions
$\tilde{\h}_\ll(X;q)$ generalise the modified Hall-Littlewood
functions $Q'(X;q)$. We shall however be concerned mainly with the
functions $\g_\ll(X;q)$.
\end{remark}

\section{Global Bases and Ribbon Functions}
\label{sec:canrib} In this section we reproduce (in slightly more
generality) a calculation due to Leclerc and Thibon \cite{LT}
which relates the $q$-Littlewood Richardson coefficients to the
global bases $G_{n\ll}$ of $\f$.  The $q$-Littlewood Richardson
coefficients $c^\ll_\mu(q)$ are defined by the expansion of the
ribbon functions on the Schur basis.

\begin{definition}
Let $\ll/\mu$ be a skew shape tileable by $n$-ribbons.  Define the
polynomials $c^{\nu}_{\ll/\mu}(q)$ by
\[
\g_{\ll/\mu}(X;q) = \sum_\nu c^{\nu}_{\ll/\mu}(q) s_\nu(X).
\]
As usual, we abuse notation by writing $c^{\nu}_{\ll}(q)$ for
$c^{\nu}_{\ll/\tll}(q)$ for partitions $\ll$ with non-empty
$n$-core.
\end{definition}
When $\ll$ has empty $n$-core, the $q$-Littlewood Richardson
coefficients are usually written in terms of the $n$-quotient:
\[
c^{\nu}_{\ll}(q) = c^{\nu}_{\ll^{(0)},\ldots,\ll^{(n-1)}}(q).
\]

We can connect the $q$-Littlewood Richardson coefficients to the
Heisenberg algebra immediately.
\begin{lem}
\label{lem:canlr} Let $\ll$ and $\mu$ be partitions.  Then
\[
S_\ll \br{\mu} = \sum_\nu c^\ll_{\nu/\mu}(-q^{-1}) \br{\nu}
\]
where the sum is over all partitions $\nu$ such that $\nu/\mu$ is
tileable by $n$-ribbons.
\end{lem}
\begin{proof}
We have by definition
\[
\g_{\nu/\mu}(X;q) = \sum_\rho \kk_{\nu/\mu,\rho}(q) m_\rho =
\sum_\ll \brac{\sum_\rho
\kk_{\nu/\mu,\rho}(q)\kappa_{\ll\rho}}s_\ll
\]
where $\kappa_{\ll\rho}$ is given by
\[
m_\rho = \sum_\ll \kappa_{\ll\rho} s_\ll.
\]
Thus
\[
c^\ll_{\nu/\mu}(q) = \sum_\rho
\kk_{{\nu/\mu},\rho}(q)\kappa_{\ll\rho}.
\] By standard results in symmetric function theory we also have
$s_\ll = \sum_\rho \kappa_{\ll\rho} h_\rho$.  Hence using
Proposition \ref{prop:hor},
\begin{align*}
S_\ll \br{\mu} &= \sum_\rho \kappa_{\ll\rho}\v_\rho \br{\mu} \\
&= \sum_\nu \brac{ \sum_\rho \kappa_{\ll\rho}
\kk_{\nu/\mu,\rho}(-q^{-1})} \br{\nu}\\
&= \sum_\nu c^\ll_{\nu/\mu}(-q^{-1}) \br{\nu}.
\end{align*}
\end{proof}

\begin{corollary}
\label{cor:canlr} Let $\ll$ be a partition.  Then
\[
G_{n\ll} = \sum_\mu c^{\ll}_\mu(-q^{-1}) \br{\mu}
\]
summed over partitions $\mu$ with no $n$-core.
\end{corollary}
\begin{proof}
This follows from Lemma \ref{lem:canlr} and Proposition
\ref{prop:can}.
\end{proof}

By Remark \ref{rem:can}, we also see that the polynomials
$c^\ll_\mu(q) = c^\ll_{\mu^{(0)},\ldots,\mu^{(n-1)}}(q)$ have
non-negative coefficients.  For $n=2$, there is a combinatorial
interpretation for the coefficients in terms of Yamanouchi domino
tableaux (see \cite{CL}).

\section{The Murnagham-Nakayama Rule}
\label{sec:mnrule}
The core calculation of this paper will be the ribbon
Murnagham-Nakayama Rule, which is essentially a consequence of the
fact that the $B_k$ act on $\f$ as a copy of the Heisenberg
algebra.  We will begin by reminding the reader of the classical
Murnagham-Nakayama Rule.

Let
\[
s_\ll(X) = \sum_\mu z_\mu^{-1}\chi_\mu^\ll p_\mu
\]
be the expansion of the Schur functions in the power sum basis.
When $\mu = (k)$ has only one part then we will write $\chi_k^\ll$
for $\chi_\mu^\ll$.  The coefficients $\chi_\mu^\ll$ are the
values of the character of $S_{|\ll|}$ indexed by $\ll$ on the
conjugacy class indexed by $\mu$.  The classical
Murnagham-Nakayama rule gives a combinatorial interpretation of
these numbers: \[ \chi_\mu^\ll = \sum_T (-1)^{h(T)} \] where the
sum is over all border-strip tableaux of shape $\ll$ and type
$\mu$.  The numbers $\chi_\mu^\ll$ are in fact the characters of
the irreducible representation labelled by $\ll$ of the symmetric
group $S_{|\ll|}$, where $\mu$ is the type of the conjugacy class.
Thus in particular all the irreducible characters of the symmetric
groups take values in the integers.  See for example \cite[Ch
7.18]{EC2}.

More generally, we have (see \cite{EC2, Mac})
\begin{prop}
\label{prop:mnrule} Let $\ll$ be a partition and $\a$ be a
composition.  Expand
\[
p_\a s_\ll(X) = \sum_\mu \chi_\a^{\mu/\ll} s_\mu(X).
\]
Then $\chi_\a^{\mu/\ll}$ is given by
\[
\chi_\a^{\mu/\ll} = \sum_T (-1)^{h(T)}
\]
where the sum is over all border strip tableaux $T$ of shape
$\mu/\ll$ and type $\a$.  Note: the border strip tableaux here
should not be confused with ribbon tableaux.  A border strip
tableaux may have border strips of different sizes.  A ribbon
tableaux has all ribbons of length $n$.
\end{prop}

This result is usually shown algebraically using the expression of
the Schur function as a bialternant $s_\ll = a_{\ll+\d}/a_{\d}$.
Theorem \ref{thm:mnviapieri} implies that it is in fact a formal
combinatorial consequence of the Pieri formula.

We will now give the analogue of Proposition \ref{prop:mnrule} for
ribbon functions.
\medskip

\begin{thm}[Ribbon Murnagham-Nakayama Rule]
\label{thm:mnrule} Let $k \geq 1$ be an integer and $\nu$ be a
partition.  Then \be{eq:mn}
\brac{1+q^{2k}+\cdots+q^{2k(n-1)}}p_k\g_\nu(X;q) = \sum_\mu
\xx^{\mu/\nu}_k(q) \g_\mu(X;q). \ee Also
\[
k \frac{\partial}{\partial p_k}\g_\nu(X;q) = \sum_\mu
\xx^{\nu/\mu}_k(q) \g_\mu(X;q).
\]
\end{thm}

\begin{proof}
Let $\d = \tilde{\nu}$ be the $n$-core of $\nu$, which we fix
throughout.  Note that the only terms $\mu$ which occur in
(\ref{eq:mn}) satisfy $\tilde{\mu} = \tilde{\nu}$.  Recall that we
will often be writing $\mu$ instead of $\mu/\d$ for convenience.

We will calculate the expression $B_{k}S_\ll\br{\d}$ with $k \geq
1$ in two ways. By Lemma \ref{lem:canlr} we can write
\[
B_{k}S_\ll\br{\d} = \sum_{\mu \in \p_\d} c^\ll_{\mu}(-q^{-1}) B_{k}\br{\mu}
\]
and by Proposition \ref{prop:bkaction} this can be written as
\[
\sum_\mu c^\ll_{\mu}(-q^{-1}) \brac{\sum_\nu
\xx^{\mu/\nu}_k(-q^{-1}) \br{\nu}} = \sum_\nu \brac{\sum_\mu
c^\ll_\mu(-q^{-1}) \xx^{\mu/\nu}_k(-q^{-1})}  \br{\nu}.
\]

On the other hand, we know by Theorem \ref{thm:heisenberg} that
$B_k$ and $S_\ll$ are both operators within a copy of the
Heisenberg Algebra.  Thus we can compute $B_k S_\ll$ within $H$.
Write
\[
S_\ll = \sum_\mu \Chi^\ll_\mu B_{-\mu}.
\]
By Lemma \ref{lem:bkcommute} and the fact that
$B_k \br{\d} = 0$ (Proposition \ref{prop:bkaction}) we have
\[
B_k S_\ll \br{\d} = \brac{\frac{1-q^{-2nk}}{1-q^{-2k}}} \sum_\mu
\Chi^{\ll/\mu}_k S_{\mu} \br{\d}
\]
since $p_k^\perp s_\ll =\sum_\mu \Chi^{\ll/\mu}_k s_\mu$ in $\Lambda$.
By Lemma \ref{lem:canlr} again, we find that this is equal to
\begin{align*}
\brac{\frac{1-q^{-2nk}}{1-q^{-2k}}} \sum_\mu \Chi^{\ll/\mu}_k
\sum_{\nu \in \p_\d} c^\mu_\nu(-q^{-1}) \br{\nu}=
 \sum_{\nu \in \p_\d} \brac{\brac{\frac{1-q^{-2nk}}{1-q^{-2k}}}
\sum_\mu \Chi^{\ll/\mu}_k c^\mu_\nu(-q^{-1})} \br{\nu}.
\end{align*}
Equating coefficients of $\br{\nu}$ we obtain \be{eq:coeff}
\brac{\frac{1-q^{-2nk}}{1-q^{-2k}}} \sum_\mu \Chi^{\ll/\mu}_k
c^\mu_\nu(-q^{-1}) = \sum_\mu c^\ll_\mu(-q^{-1})
\xx^{\mu/\nu}_k(-q^{-1}). \ee

We now calculate
\begin{align*}
\brac{\frac{1-q^{2nk}}{1-q^{2k}}}p_k\g_\nu(X;q) &= \brac{\frac{1-q^{2nk}}{1-q^{2k}}} \sum_\mu c^\mu_\nu(q) p_k s_\mu \\
&= \brac{\frac{1-q^{2nk}}{1-q^{2k}}} \sum_\mu c^\mu_\nu(q)
\brac{\sum_\ll \Chi^{\ll/\mu}_k s_\ll} \\
&= \sum_\ll \brac{\brac{\frac{1-q^{2nk}}{1-q^{2k}}} \sum_\mu
c^\mu_\nu(q) \Chi^{\ll/\mu}_k} s_\ll \\
& = \sum_\ll \brac{ \sum_{\mu \in \p_\d} c^\ll_\mu(q)
\xx^{\mu/\nu}_k(q)}
s_\ll \;\;\; \mbox{ using Equation (\ref{eq:coeff})} \\
&= \sum_{\mu \in \p_\d} \xx^{\mu/\nu}_k(q) \g_\mu(X;q).
\end{align*}

This proves the first statement.  The second statement is proved
in the same manner, considering $B_{-k}$ instead of $B_k$.
\end{proof}
We shall see later in Section \ref{sec:skew} that the lowering
version of the Murnagham-Nakayama rule can be deduced
combinatorially in a rather straightforward manner (using an
observation of Schilling, Shimozono and White \cite{SSW}). In fact
it is clear that the lowering operator version is easier as the
proof does not require the use of the commutator relation of $B_k$
and $B_{-k}$ in Theorem \ref{thm:heisenberg}, only that the
$B_{-k}$ for $k \geq 1$ commute.  Thus the lowering version of
Theorem \ref{thm:mnrule} is essentially equivalent to the fact
that the $\g_\ll(X;q)$ are symmetric.

Note that it is rather difficult to interpret Theorem
\ref{thm:mnrule} in terms of the $n$-quotient at $q=1$.  When
$q=1$ the product
$\brac{1+q^{2k}+\cdots+q^{2k(n-1)}}p_k\g_\ll(X;q)$ becomes
\[
np_k s_{\ll^{(0)}}s_{\ll^{(1)}}\cdots s_{\ll^{(n-1)}}
\]
which may be written as the sum of $n$ usual Murnagham-Nakayama
rules as
\[
\sum_{i=0}^{n-1} s_{\ll^{(0)}}\cdots (p_k s_{\ll^{(i)}}) \cdots
s_{\ll^{(n-1)}}.
\]
Thus we might expect that border ribbon strips of size $k$
correspond to adding a usual ribbon strip of size $k$ to one
partition in the $n$-quotient.  However, the following example
shows that this will not work.

\begin{example}
Let $n = 2$ and consider $(1+q^4)p_2 \cdot 1$.  By the ribbon
Murnagham-Nakyama rule ($\g_0 = 1$), this should equal to
\[
\g_{(4)}+q\g_{(3,1)}+(q^2-1)\g_{(2,2)}-q\g_{(2,1,1)}-q^2\g_{(1,1,1,1)}.
\]
We can compute directly that
\begin{align*}
\g_{(4)} = h_2,\;\;\;\g_{(3,1)} =qh_2, \;\;\; \g_{(2,1,1)} = qe_2 \\
\g_{(2,2)} = q^2h_2+e_2, \;\;\; \g_{(1,1,1,1)} = q^2 e_2,
\end{align*}
verifying Theorem \ref{thm:mnrule} directly.  On the other hand,
the shapes which correspond to a single border strip in one partition
of the 2-quotient are
$\set{(4),(3,1),(2,1,1),(1,1,1,1)}$ and the corresponding
$\g_\ll$ terms do not give $(1+q^4)p_2$.
\end{example}
It seems possible that the ribbon Murnagham-Nakayama rule may
have some relationship with the representation theory of the
wreath products $S_n\S C_p$, or even more likely to the cyclotomic
Hecke algebras associated to these wreath products (see for
example \cite{Mat}).

\section{The map $\m: \f \rightarrow \Lambda(q)$}
\label{sec:map} In this section we will study a linear map from
$\f$ to $\Lambda(q)$.

\begin{definition}
Let $\m: \f \rightarrow \Lambda(q)$ be the linear over
$\C$ map defined by $q \mapsto -q^{-1}$ and
\[
\br{\ll} \mapsto \g_\ll.
\]
\end{definition}

As the $\br{\ll}$ form a basis of $\f$ it is clear that such a map
exists and is unique.  Since the $\g_\ll$ span $\Lambda(q)$ but
are not linearly dependent, this map is surjective but not
injective.  We would intuitively think that in doing so we have
lost a lot of information by going from $\f$ to $\Lambda(q)$ but
curiously this map has many remarkable properties.  Note that the
map $\m$ should not be confused with the classical identification
of $\f$ with $\Lambda$ via $\br{\ll} \leftrightarrow s_\ll$, which we shall
comment about in Section \ref{sec:open}.

The following theorem can be interpreted as saying that $\m$ is a projection
of $\f$ onto $\CH$.
\begin{thm}
\label{thm:map}
After changing $q$ to $-q^{-1}$, the map $\m$ is a map of $H$ modules.  More precisely,
let $A$ be an element of the Heisenberg algebra
$H$.  Identify
$\m(A)$ with the element of $\Lambda(q)$ given by
\begin{align*}
B_k &\mapsto k \frac{\partial}{\partial p_k} \\
B_{-k} &\mapsto \brac{\frac{1-q^{2nk}}{1-q^{2k}}}p_k. \\
\end{align*}
(Recall from Section \ref{sec:hei} that this identifies $\Lambda(q)$
with $\CH$ and defines an action of $H$ with $q \mapsto -q^{-1}$.)
Then we have
\[
\m(A\cdot v) = \m(A)\m(v)
\]
for any $v \in \f$.

Furthermore, we have \[\m(G_{n\ll}) = s_\ll[(1+q^2+\cdots
+q^{2n-2})X].\]
\end{thm}
\begin{proof}
The first claim follows from Theorem \ref{thm:mnrule} and
Proposition \ref{prop:bkaction} as we can simply compare both
sides of the equation on the spanning sets $\br{\ll}$ and
$\g_\ll$.  To obtain the second claim, we apply the first claim
with $v = \br{0}$ and $A = S_\ll$, using also $G_{n\ll}=S_\ll
\br{0}$ by Proposition \ref{prop:can}.
\end{proof}

In later sections we shall see that the two involutions of $\f$,
the bar involution and the involution $v \mapsto v'$ become
algebra isomorphisms of $\CH \cong \Lambda(q)$.  While most of the
later results can be phrased concisely in terms of the
Heisenberg algebra, we shall continue to use symmetric function
terminology, thinking of the $\g_\ll(X;q)$ as elements in
$\Lambda(q)$ rather than $\CH$.

\medskip

We also have the following explicit descriptions of
$s_\ll[(1+q^2+\cdots +q^{2n-2})X]$.

\begin{corollary}
\label{cor:phican} Let $\ll$ be a partition and $\d$ a fixed
$n$-core. In $\Lambda(q)$ we have
\[
s_\ll[(1+q^2+\cdots +q^{2n-2})X] = \sum_\mu c^\ll_\mu(q)
\g_\mu(X;q)
\]
and
\[
s_\ll[(1+q^2+\cdots +q^{2n-2})X] = \sum_{\mu,\nu} c^\ll_\mu(q)
c^\nu_\mu(q) s_\nu(X)
\]
where the sums are over all partitions $\mu$ with $n$-core $\d$.
\end{corollary}
\begin{proof}
These are immediate consequences of Theorem \ref{thm:map} and
Lemma \ref{lem:canlr} as $\m(\br{\d}) = 1$ for an $n$-core $\d$.
\end{proof}
By Corollary \ref{cor:canlr} we have calculated the images of
those global basis vectors which are highest weight vectors for
$\uqsln$ in $\f$. It is not clear whether this leads to any
interesting results concerning $G_{n\ll}$ in $\f$.

Applying $\m$ to both sides of Lemma \ref{lem:canlr} we obtain
\begin{prop}
\label{prop:canlr}
Let $\ll$ and $\mu$ be partitions. Then
\[
s_\ll[(1+q^2+\cdots +q^{2n-2})X]\g_\mu(X;q)= \sum_\nu c^{\ll}_{\nu/\mu}(q)\g_\nu(X;q).
\]
\end{prop}
This could be thought of as some kind of Littlewood Richardson rule for
ribbon tableaux.  In fact the coefficients $c^{\ll}_{\nu/\mu}(q)$ are
the $q$-Littlewood Richardson coefficients which are the coefficients
of the expansion of $\g_{\nu/\mu}(X;q)$ on the Schur basis.  Of course
there is no combinatorial description of these coefficients except
in the case $n=2$, via the Yamanouchi domino tableaux of Carr\'{e} and Leclerc \cite{CL}.

\section{Ribbon Pieri Formulae}
\label{sec:pieri} Recall that in the classical theory of symmetric
functions and in the enumerative geometry of the Grassmanian we
have the formula
\[
h_k s_\ll = \sum_\mu s_\mu
\]
where the sum is over all $\mu$ such that $\mu/\ll$ is a
horizontal strip with $k$ boxes.  This formula arises in the
intersection theory of the Grassmanian as the intersection of a
special Schubert class with an arbitrary Schubert class and was
discovered geometrically by Pieri \cite{Pie}.  We will now give
the ribbon analogue of the Pieri formula.

\medskip

Let $n \geq 1$ be a fixed integer.  Define the formal power series
\[
H(t) = \prod_i \prod_{k=0}^{n-1} \frac{1}{1-x_i q^{2k} t}
\]
and
\[
E(t) = \prod_i \prod_{k=0}^{n-1} \brac{1+x_i q^{2k} t}.
\]
As usual we may define
symmetric functions $\hh_k$ and $\e_k$ by
\[
H(t) = \sum_k \hh_k t^k
\]
and
\[
E(t) = \sum_k \e_k t^k.
\]
Note that we have suppressed the integer $n$ from the notation. We
shall see later that the definitions of these power series are
completely natural in the context of Robinson-Schensted ribbon
insertion.

\medskip
In plethystic notation, $\hh_k = h_k[(1+ q^{2} +
\cdots+q^{(2n-2)})X]$ and $\e_k =e_k[(1+ q^{2} +
\cdots+q^{(2n-2)})X]$.  This can be seen as follows.  Write
\begin{align*}
\log(\sum_k h_k t^k) &= \sum_i \log \brac{\frac{1}{1-x_i t}} \\
&=  \sum_i \sum_r \frac{(x_it)^r}{r} \\
&= \sum_r \frac{p_r t^r}{r}. \\
\end{align*}
Now take the plethysm $\plqn(p_r) = (1+ q^{2r} +
\cdots+q^{(2n-2)r})p_r$ and reverse all the steps.
\medskip

The following theorem is an immediate consequence of Theorem
\ref{thm:mnviapieri} and Theorem \ref{thm:mnrule}.  Alternatively,
one could just use Theorem \ref{thm:map} and Proposition
\ref{prop:hor}.

\begin{thm}[Ribbon Pieri Rule]
\label{thm:pieri} Let $\lambda$ be a partition.  Then
\be{eq:pieri} \hh_k \g_\lambda(X;q) = \sum_\mu q^{s(\mu/\lambda)}
\g_\mu(X;q) \ee where the sum is over all partitions $\mu$ such
that $\mu / \lambda$ is a horizontal $n$-ribbon strip with $k$
ribbons. Here $s(\mu/\lambda)$ refers to the spin of the unique
tableaux which is a horizontal ribbon strip of shape
$\mu/\lambda$. Also
\[
\e_k \g_\lambda(X;q) = \sum_\mu q^{s(\mu/\lambda)} \g_\mu(X;q)
\]
where the sum is over all partitions $\mu$ such that $\mu /
\lambda$ is a vertical $n$-ribbon strip with $k$ ribbons. Here
$s(\mu/\lambda)$ refers to the spin of the unique tableaux which
is a vertical ribbon strip of shape $\mu/\lambda$.
\end{thm}
One can obviously obtain the corresponding formulae for
$\hh_\a = \hh_{\a_1}\cdots\hh_{\a_r}$ in terms of ribbon tableaux
with weight $\a$.

We can also obtain the two statements of Theorem \ref{thm:pieri}
from each other via the involution $\w_n$ of Section \ref{sec:inv}.
Note that by Theorem \ref{thm:pieri}, we have
\[
\hh_k = \sum_\ll q^{\text{mspin}(\ll)}\g_\ll(X;q)
\]
where the sum is over all $\ll$ with no $n$-core such that $|\ll|
= kn$ with no more than $n$ rows.  Applying the usual Cauchy identity 
one sees that
\[
\hh_k = \sum_{|\mu|=k} s_\mu(1,q^2,\ldots,q^{2(n-1)})s_\mu(X).
\]
Taking the coefficient of $p_1^n$ on both sides we see that a modified spin
generating function of ribbon tableaux $T$ of size $k$ and shape $\ll$
satisfying $\tll = \emptyset$ and $l(\ll) \leq n$ is
\[
\sum_T q^{\text{mspin}(sh(T))}q^{s(T)} = \sum_{|\mu|=k}s_\mu(1,q^2,\ldots,q^{2(n-1)})f^{\mu}
\]
where $f^{\mu}$ denotes the number of standard Young tableaux of shape $\mu$.

\begin{example}
Let $n=3$, $k=2$ and $\ll = (3,1)$.  Then
\[
\hh_2 \g_{(3,1)} = \g_{(9,1)} +q \g_{(6,2,2)} + q^2\g_{(4,4,2)} +
q^2\g_{(6,1,1,1,1)} + q^3 \g_{(3,3,2,1,1)} + q^4 \g_{(3,2,2,2,1)}.
\]
\end{example}
\medskip

Setting $q=1$ in $H(t)$ we see that
\[
\hh_k(X;1) = \sum_\a h_\a
\]
where the sum is over all compositions $\a =
(\a_0,\ldots,\a_{n-1})$ satisfying $\a_0 + \cdots + \a_{n-1} = k$.
We may thus interpret Theorem \ref{thm:pieri} at $q=1$ in terms of
the $n$-quotient as the following formula: \be{eq:pieriq1}
\brac{\sum_\a h_\a} s_{\ll^{(0)}}\cdots s_{\ll^{(n-1)}} = \sum_\a
\brac{h_{\a_0}s_{\ll^{(0)}}}\cdots
\brac{h_{\a_{n-1}}s_{\ll^{(n-1)}}} \ee where the sum is over the
same set of compositions as above. Note that the right hand side
of (\ref{eq:pieriq1}) is indeed equal to the right hand side of
(\ref{eq:pieri}) at $q=1$ since a horizontal ribbon strip of size
$k$ is just a union of horizontal strips with total size $k$ in
the $n$-quotient.

It is clear that we also obtain lowering versions of the Pieri
rules.  If $h_k = f(p_1, p_2, \ldots)$ we know that the adjoint
operator (with respect to the usual inner product) is $h_k^\perp =
f(\diff{p_1}, 2\diff{p_2}, \ldots)$.  Thus by Theorem
\ref{thm:mnrule} and Theorem \ref{thm:mnviapieri} we have
\begin{prop}[Ribbon Pieri Rule -- Lowering Version]
\label{prop:lowpieri}
Let $\ll$ be a partition and $k \geq 1$ be an integer.  Then
\[
h_k^\perp \g_\ll(X;q) = \sum_\mu q^{s(\ll/\mu)} \g_\mu(X;q)
\]
where the sum is over all $\mu$ such that $\ll/\mu$ is a
horizontal ribbon strip and $s(\ll/\mu)$ is the spin of such a
horizontal ribbon strip.  Similarly,
\[
e_k^\perp \g_\ll(X;q) = \sum_\mu q^{s(\ll/\mu)} \g_\mu(X;q)
\]
where the sum is over all $\mu$ such that $\ll/\mu$ is a vertical
ribbon strip and $s(\ll/\mu)$ is the spin of such a vertical
ribbon strip.
\end{prop}
This is a spin version of a branching formula first observed by
Schilling, Shimozono and White \cite{SSW} (see Section
\ref{sec:skew}).

\section{The Ribbon Involution $\w_n$}
\label{sec:inv}

In this section we will define an involution $w_n$ on $\Lambda(q)$
which is essentially the involution $v \mapsto v'$ on the Fock
space $\f$ of Section \ref{sec:can}.  However, this involution
will turn out to be not just a semi-linear involution, but also a
$\C$-algebra isomorphism of $\Lambda(q)$.

\begin{definition}
Define the \emph{ribbon involution} $w_n: \Lambda(q) \rightarrow
\Lambda(q)$ as the semi-linear map satisfying $w_n(q) = q^{-1}$
and
\[
w_n(s_\lambda) = q^{(n-1)|\lambda|}s_{\lambda'}.\]
\end{definition}

\begin{thm}
\label{thm:inv} The map $w_n$ is an $\C$-algebra homomorphism
which is an involution.  It maps $\g_{\ll/\mu}$ into
$\g_{(\ll/\mu)'}$ for every skew shape $\ll/\mu$.
\end{thm}
\begin{proof}
The fact that $w_n$ is an algebra homomorphism follows from the
fact that if $s_\ll s_\mu = \sum c_{\ll \mu}^{\nu} s_\nu$ then
$s_{\ll'} s_{\mu'} = \sum c_{\ll \mu}^{\nu} s_{\nu'}$, and that
the grading is preserved by multiplication. That $w_n$ is an
involution is a quick calculation.

For the last statement, we use Proposition \ref{prop:semilinear}
and Lemma \ref{lem:canlr} which give
\begin{align*}
\brac{S_\nu\br{\mu}}' &= (-q)^{(n-1)k}S_{\nu'}\br{\mu'} \\
\sum_\ll c^\nu_{\ll/\mu}(-q) \br{\ll'} & = (-q)^{(n-1)k}\sum_\ll
c^{\nu'}_{\ll'/\mu'}(-q^{-1}) \br{\ll'}.
\end{align*}
Here $k = |\nu|$.
Equating coefficients of $\br{\ll'}$ and changing $q$ to $-q^{-1}$ we obtain
\[
c^\nu_{\ll/\mu}(q^{-1}) = q^{-(n-1)k} c^{\nu'}_{\ll'/\mu'}(q).
\]
Thus
\begin{align*}
w_n(\g_{\ll/\mu}) &= \sum_\nu w_n(c_{\ll/\mu}^\nu(q)s_\nu) \\
                &= \sum_\nu \brac{c^{\nu'}_{\ll'/\mu'}(q)q^{-(n-1)|\nu|}}q^{(n-1)|\nu|}s_{\nu'}\\
                 &= \g_{\ll'/\mu'}.
\end{align*}
\end{proof}

\begin{prop}
\label{prop:wncommute} Let $f \in \Lambda(q)$ have degree $k$.
Then we have
\[
q^{2(n-1)k}\w_n \brac{\plqn(f)} = \plqn\brac{\w_n(f)}.
\]
In particular,
\[
\w_n \brac{s_\ll[(1+q^2+\cdots+q^{2(n-1)})X]} =
q^{-(n-1)k}s_\ll[(1+q^2+\cdots+q^{2(n-1)})X].
\]
\end{prop}
\begin{proof}
Since both $\w_p$ and $\plqn(f)$ are $\C$-algebra homomorphisms we
need only check this for the elements $p_k$ and for $q$, for which
the computation is straightforward.
\end{proof}

\section{The Ribbon Cauchy Identity}
\label{sec:cauchy} It is well known (see \cite{Mac,EC2}) that the
Schur functions satisfy the equation \be{eq:classCauchy} \sum_{\ll
\in \p} s_{\ll}(X)s_{\ll}(Y) = \prod_{i,j} \frac{1}{1-x_iy_j} \ee
known as the Cauchy identity.  This is equivalent to the fact that
the Schur functions form an orthonormal basis of $\Lambda$.  It
also describes the decomposition of the polynomial functions on
the $m \times m$ matrices under the (commuting) left and right
actions of $GL_m$. Equation (\ref{eq:classCauchy}) is often proven
combinatorially via the Knuth's extension of the
Robinson-Schensted correspondence, which is a bijection between
matrices with certain row and column sums and pairs of
semistandard Young tableaux of the same shape. This combinatorial
approach to the Cauchy identity for ribbon tableaux will be
studied in Section \ref{sec:ribinsertion}.

Let us write the formal power series
\[
\o(X;q) = \prod_{i,j}  \prod_{k=0}^{n-1} \frac{1}{1-x_i y_j q^{2k}}.
\]
A dual version of this series is
\[
\to(X;q) = \prod_{i,j}  \prod_{k=0}^{n-1} \brac{1+x_i y_j q^{2k}}.
\]

\begin{thm}[Ribbon Cauchy Identity]
\label{thm:cauchy} Fix $n$ as usual and a $n$-core $\d$.  Then
\[
\o(X;q) = \sum \g_\lambda(X;q) \g_\lambda(Y;q)
\]
where the sum is over all $\ll$ such that $\tilde{\lambda} = \d$.
\end{thm}
Note that this does not imply that the $\g_\ll$ form an
orthonormal basis under a certain inner product, as they are not
linearly independent.
\begin{proof}
By Corollary \ref{cor:phican} we have
\[
s_\ll[(1+q^2+\cdots +q^{2n-2})X] = \sum_\mu c^\ll_\mu(q)
\g_\mu(X;q)
\]
where the sum is over all $\mu \in \p_\d$.  Thus
\begin{align*}
\sum_\ll s_\ll[(1+q^2+\cdots +q^{2n-2})X] s_\ll(Y) &= \sum_\mu
\brac{\sum_\ll c^\ll_\mu(q) s_\ll(Y)} \g_\mu(X;q) \\
&= \sum_\mu \g_\mu(X;q) \g_\mu(Y;q).
\end{align*}

Let $\plqn(X)$ denote the algebra automorphism of $\Lambda[X](q)
\otimes_{\C(q)} \Lambda[Y](q)$ given by \[p_k(X) \mapsto
(1+q^{2k}+\cdots+q^{(2n-2)k})p_k(X).\]

Applying $\plqn(X)$ to
\[
\log\brac{\prod_{i,j} \frac{1}{1-x_iy_j}} = \sum_k
\frac{1}{n}p_k(X)p_k(Y)
\]
gives
\[
\log\brac{\prod_{i,j} \prod_{k=1}^{n-1} \frac{1}{1-x_iy_jq^{2k}}}
\]
which is exactly $\log(\o)$.  But applying $\plqn(X)$ to the left
hand side of (\ref{eq:classCauchy}) gives
\[
\sum_\ll s_\ll[(1+q^2+\cdots +q^{2n-2})X] s_\ll(Y)
\]
from which the Theorem follows.
\end{proof}

Now let us compute $\w_n(\o)$ where we let
\[
\w_n: \Lambda[X](q) \otimes_{\C(q)} \Lambda[Y](q) \rightarrow
\Lambda[X](q) \otimes_{\C(q)} \Lambda[Y](q)\] act on the $X$
variables by
\[
\w_n(f(X;q) \otimes g(Y;q)) \mapsto \w_n(f(X;q)) \otimes
g(Y;q^{-1}).
\]
One checks immediately that this is indeed an algebra involution.
We have (fixing an $n$-core $\d$)
\begin{align*}
\w_n(\o) &= \w_n\brac{\sum_{\ll \in \p_\d} \g_\ll(X;q)\g_\ll(Y;q)} \\
&= \sum_{\ll \in \p_\d} \g_{\ll'}(X;q)\g_\ll(Y;q^{-1}).
\end{align*}

Also,
\begin{align*}
\w_n(\o) & = \w_n\brac{\sum_\ll s_\ll(X)
s_\ll[(1+q^2+\cdots+q^{2(n-1)})Y]} \\
&= \sum_\ll q^{(n-1)|\ll|}s_{\ll'}(X)s_\ll[(1+q^{-2}+\cdots+q^{-2(n-1)})Y] \\
&= \prod_{i,j} \prod_{k=0}^{n-1} (1+x_i y_j q^{n-1-2k}).\\
\end{align*}

Thus
\[
\sum_{\ll \in \p_\d} \g_{\ll'}(X;q)\g_\ll(Y;q^{-1}) = \prod_{i,j}
\prod_{k=0}^{n-1} (1+x_i y_j q^{n-1-2k}).
\]
If we multiply the $d^{th}$ graded piece of each side by
$q^{(n-1)d}$ we obtain the following result.

\begin{prop}
\label{prop:dualcauchy} Fix an $n$-core $\d$.  We have
\[
\to = \sum_{\ll \in \p_\d}
q^{(n-1)|\ll/\tll|}\g_{\ll'}(X;q)\g_\ll(Y;q^{-1}).
\]
\end{prop}
The factor of $q^{(n-1)|\ll/\tll|}$ can be explained
combinatorially by the fact that $s(T') = q^{(n-1)|\ll/\tll|} -
s(T)$ for a ribbon tableaux $T$ and its conjugate $T'$ satisfying
$sh(T) = \ll$.

\section{The Ribbon Inner Product and the Bar Involution on $\Lambda(q)$}
\label{sec:ip}
In this section we will define an inner product on $\Lambda(q)$ which seems particularly adapted
to the study of ribbon functions.  We will also give a $\C$-algebra involution on
$\Lambda(q)$ which is compatible with the bar involution of $\f$ (for at least
the space of highest weight vectors).

\begin{definition}
Let $\ip{., .}_n: \Lambda(q) \times \Lambda(q) \rightarrow \C(q)$
be the $\C(q)$-bilinear map defined by
\[
\ip{\plqn(p_{\ll}),p_\mu} = z_\ll \d_{\ll\mu}.
\]
\end{definition}
It is clear that $\ip{., .}_n$ is non-degenerate.

The inner product $\ip{\,.,\,.}_n$ is related to $\o$ in the same way as the usual
inner product is related to the usual Cauchy kernel:
\begin{prop}
Two bases $\set{v_\ll}$ and $\set{w_\ll}$ of $\Lambda(q)$ are dual with respect to
$\ip{\,.\,}_n$ if and only if
\[
\sum_\ll v_\ll(X)w_\ll(Y) = \o.
\]

In particular, $\set{s_\ll[(1+q^2+\ldots+q^{2(n-1)})X]}$ is dual to $\set{s_\ll}$.
\end{prop}
\begin{proof}
It is clear that $\set{p_{\ll}/z_\ll}$ and $\set{\plqn{p_\ll}}$ are dual.  But applying
$\plqn$ to the usual Cauchy kernel gives
\[
\sum_{\ll} \frac{1}{z_\ll}\plqn(p_\ll(X))p_\ll(Y) = \o.
\]
To see that this is true for any pair of dual bases of
$\Lambda(q)$ with respect to $\ip{\,.,\,.}_n$ is an exercise in
linear algebra (see for example \cite[Lemma 7.9.2]{EC2}). The last
statement is a consequence of Theorem \ref{thm:cauchy}.
\end{proof}
In fact if $\set{v_\ll}$ and $\set{w_\ll}$ are dual basis of $\Lambda$ with respect
to the usual inner product then it is clear that $\set{\plqn(v_\ll)}$ and $\set{w_\ll}$
are dual with respect to $\ip{\,.,\,.}_n$. We now give some basic properties of $\ip{\,., \,.}_n$.

\begin{lem}
The inner product $\ip{\,., \,.}_n$ is symmetric.
\end{lem}
\begin{proof}
This is clear from the definition as we can just check this on the basis $p_\ll$ of $\Lambda(q)$.
\end{proof}

Recall that for $f\in \Lambda$, $f^\perp$ denotes its adjoint with respect
to the usual inner product.

\begin{prop}
\label{prop:perp}
The operator $f^\perp$ is adjoint to multiplication by $\plqn(f) \in \Lambda(q)$.
\end{prop}
\begin{proof}
This is a consequence of $\ip{f,g}$ = $\ip{\plqn(f),g}_n$.
\end{proof}
The inner product $\ip{\,., \,.}_n$ is compatible with the inner
product $\ip{\br{\ll},\br{\mu}} =\d_{\ll\mu}$ on $\f$ when we
restrict our attention to the space of highest weight vectors.  In
$\f$ we have $\ip{B_ku,v} = \ip{u,B_{-k}v}$ for any $u,v \in \f$, see
\cite[Proposition 7.9]{LT} which corresponds to Proposition \ref{prop:perp}.

The bar involution ${}^-: \f \rightarrow \f$ of Section \ref{sec:can} also has an image
under $\m$.

\begin{definition}
Define the $\C$-algebra involution ${}^-: \Lambda(q) \rightarrow \Lambda(q)$ by
$\ov{q} = q^{-1}$ and
\[
p_k \mapsto q^{2(n-1)k}p_k.
\]
\end{definition}
It is clear that ${}^-$ is indeed an involution.  The following proposition shows
in particular that
${}^-: \Lambda(q) \rightarrow \Lambda(q)$ is the image of the bar involution on $\f$ under
$\m$.  This implies that $\ip{\m(u),\m(v)}_n = \ip{u,v}$ for two vectors
$u,v \in \f$ lying within the subspace of highest weight vectors.

\begin{prop}
\label{prop:barinv} Let $u,v \in \Lambda(q)$.  The involution
${}^-: \Lambda(q) \rightarrow \Lambda(q)$ has the following
properties:
\begin{align*}
\ov{\m(G_{n\ll})} &= \m(G_{n\ll}), \\
\ov{\plqn(p_k)} &= \plqn(p_k),\\
\ip{\ov{u},v}_n &= \ip{\w_n(u),\ov{\w_n(v)}}.
\end{align*}

\end{prop}
\begin{proof}
As ${}^-$ is an algebra homomorphism, the first statement follows
from the second statement and Theorem \ref{thm:map}.  The second
statement is a straightforward computation.  For the last
statement, we compute explicitly both sides for the basis $p_\ll$
of $\Lambda(q)$.
\end{proof}

Proposition \ref{prop:barinv} shows that $\ov{\m(v)} =\m(\ov{v})$
for all $u,v$ in the subspace of highest weight vectors in $\f$.  However
this is not true in general.  For example, $\br{(3,1)}+q\br{(2,2)}+q^2\br{(2,1,1)}$
is bar invariant in $\f$ but its image under $\m$ is not.

\section{Skew and super ribbon functions}
\label{sec:skew} We now describe some properties of the skew
ribbon functions $\g_{\ll/\mu}(X;q)$. Unfortunately, we have been
unable to describe them in terms of an adjoint in analogy with the
formula
\[
s_{\ll/\mu} = s_\ll^\perp s_\mu.
\]
However, the following proposition is an analogue of
\[
\sum_\ll s_\ll(X)s_{\ll/\mu}(Y) = s_\mu(X)\prod_{i,j}\frac{1}{1-x_iy_j}.
\]

\begin{prop}
Let $\mu$ be any partition.  Then
\[
\g_\mu(X;q)\o = \sum_\ll \g_\ll(X;q)\g_{\ll/\mu}(Y;q)
\]
where the sum is over all $\ll$ satisfying $\tll = \tilde{\mu}$.
\end{prop}
\begin{proof}
Lemma \ref{lem:canlr} implies that
\[
\plqn(s_\nu(X))\g_\mu(X;q) = \sum_\ll c^\nu_{\ll/\mu}(q)
\g_\ll(X;q).
\]
Now multiply both sides by $s_\nu(Y)$ and sum over $\nu$.  Finally
use Theorem \ref{thm:cauchy}.
\end{proof}

Another essentially equivalent way in which skew ribbon functions
arise was observed by Schilling, Shimozono and White \cite{SSW} in
cospin form.  By the combinatorial definition of $\g_\ll$ we immediately have the
coproduct expansion
\[
\g_\ll(X+Y;q) = \sum_\mu \g_\mu(X;q)\g_{\ll/\mu}(Y;q).
\]
Since (\cite{Mac})
\[
\Delta f = \sum_\mu s_\mu^\perp f \otimes s_\mu
\]
we get immediately that
\[
s_\nu^\perp \g_\ll(X;q) = \sum_\mu \g_\mu(X;q)
\ip{\g_{\ll/\mu}(Y;q),s_\nu}.
\]
Setting $\nu = (k)$ we obtain the lowering version of the Pieri rule (Proposition
\ref{prop:lowpieri}):
\[
h_k^\perp \g_\ll(X;q) \sum_\mu q^{s(\ll/\mu)}\g_\mu(X;q)
\]
where the sum is over all $\mu$ such that $\ll/\mu$ is a
horizontal ribbon strip of size $k$.
\medskip

We would like to mention another generalisation of the usual ribbon functions 
which are super ribbon functions.  Fix a total order on two alphabets $A =
\set{1 < 2 < 3 < \cdots}$ and $A' =\set{1' < 2' < 3' < \cdots}$
(which we assume to be compatible with each of their natural
orders).

\begin{definition}
A super ribbon tableaux $T$ of shape $\ll/\mu$ is a ribbon
tableaux of the same shape with ribbons labelled by the two
alphabets such that the skew shape containing ribbons labelled $a
\in A$ form a horizontal ribbon strip and those labelled $a' \in
A'$ form a vertical ribbon strip.  These strips are required to be
compatible with the chosen total order on $A \cup A'$, as usual.

Define the \emph{super ribbon
function} $\g_{\ll/\mu}(X/Y;q)$ as the following weight and spin 
generating function:
\[
\g_{\ll/\mu}(X/Y;q) = \sum_T q^{s(T)}{\mathbf x}^{w(T)}{\mathbf
y}^{w'(T)}
\]
where the sum is over all super ribbon tableaux $T$ of shape
$\ll/\mu$ and $w(T)$ is the weight in the first alphabet $A$ while
$w'(T)$ is the weight in the second alphabet $A'$.
\end{definition}

\begin{prop}
\label{prop:srf} The super ribbon function $\g_{\ll/\mu}(X/Y;q)$
is a symmetric function in the $X$ and $Y$ variables, separately.
\end{prop}
\begin{proof}
The proof is completely analogous to that of Theorem
\ref{thm:sym}, using the commutativity of both the operators $\v_k$ and $\tv_k$.
\end{proof}
No doubt the super ribbon functions can be studied in the same way that super Schur functions are.

\section{Open questions and other aspects of ribbon functions}
\label{sec:open} In this section we describe some other aspects
of ribbon functions which we have not mentioned or which may be
interesting for further study.  The problem of finding combinatorial
proofs of the theorems in this part will be addressed in Section
\ref{sec:ribinsertion}.

\textbf{Cospin vs. spin.}
Recall that $\text{cosp}(T) = \text{mspin}(T) - s(T)$ for a ribbon
tableaux $T$.  It is easy to see that $\text{cosp}(T)$ is always even.
In many situations it appears that the statistic cospin is more natural
than the statistic spin.  For example, Lascoux, Leclerc and Thibon \cite{LLT}
have shown that the cospin $\tilde{\h}(X;q)$ functions are
generalisations of Hall-Littlewood Functions.  Cospin also appears to be
the natural statistic when finding connections between ribbon tableaux
and rigged configurations (see for example \cite{Sch}).

All the formulae in this Part can be phrased in terms
of cospin if suitable powers of $q$ are inserted.  However, it is clear
that the formulae presented in terms of spin is the more natural form.

\textbf{Vertex Operators.}
The relationship between the ribbon functions and the Heisenberg algebra
suggests that we may want to study the affect of the `vertex operators'
\[
A_r = \sum_i \hh_{r-i}\e_i^\perp
\]
on the ribbon functions.  It is well known that both the Schur functions
and the Hall Littlewood functions \cite{Jin} can be developed in the context of
these vertex operators so perhaps the ribbon functions can be studied
in the same way.

\textbf{Generalised Kostka polynomials.} There is a mysterious and
unsolved connection between the generalised Kostka polynomials of
\cite{KS, SW3, SchW} and the $q$-Littlewood Richardson
coefficients. The coefficients of the generalised Kostka
polynomials are not always positive but when the $n$-quotient is a
sequence of rectangles the $q$-Littlewood Richardson coefficients
appear to coincide with the generalised Kostka polynomials. A
vertex operator description of the generalised Kostka polynomials
has also been given by Shimozono and Zabrocki \cite{SZ}.  Perhaps
the vertex operators described above will help.

\textbf{Jacobi-Trudi and alternant formulae.}
The Pieri rule (Theorem \ref{thm:pieri}) we have given stops short of
giving a closed formula for the functions $\g_\ll(X;q)$.  It is
well known that the Schur functions can be written as
\[
s_\ll = \frac{a_{\ll+\d}}{a_\d}
\]
and as
\[
s_\ll = det(h_{\ll_i - i + j})_{i,j=1}^{l(\ll)}.
\]
Most algebraic treatments (see \cite{Mac}) of the theory of
symmetric functions use these formulae as the basis of all the
algebraic computations for Schur functions.  It would be nice to
have a similar closed formula for the ribbon functions.

\textbf{Enumerative problems.} Stanley \cite{EC2} has given a
`hook content formula' for the specialisation
$s_\ll(1,t,t^2,\ldots,t^r)$ of the Schur functions. In particular
this gives the hook length formula for the number of standard
Young tableaux of a particular shape.  At $q=1$ the corresponding
problem for ribbon tableaux is trivial due to Littlewood's
$n$-quotient map.  However, can anything be done for arbitrary
$q$?

When $n=2$, the specialisation $q^2 = -1$ relates domino tableaux
to the study of enumerative study of sign-imbalance
\cite{Sta,Whi,Lam}. It is not clear whether this can be
generalised to arbitrary $n$.

\textbf{Graded $S_n$ representations and $\hh_k$.}  The
non-negativity of the $q$-Littlewood Richardson coefficients
$c^\mu_\ll(q)$ would follow from the existence of a graded $S_n$
representation with Frobenius character $\g_\ll(X;q)$ (where the
coefficient of powers of $q$ correspond to the graded parts).

Such a graded $S_n$ representation can be easily described for the
ribbon homogeneous function $\hh_k$ (see \cite[Ex. 7.75]{EC2}).
Let $S_k$ act on the multiset $M =
\set{1^{n-1},2^{n-1},\ldots,k^{n-1}}$ in the natural way.  Then
the representation we seek is given by the action of $S_k$ on the
subsets of $M$ with the grading given by the size of such a
subset.  This suggests that one might seek subrepresentations of
this representation which correspond to the $\g_\ll(X;q)$ for $l(\ll) \leq n$.

\textbf{Connections with diagonal harmonics.}  When $k = n-1$, the
function $\e_k$ also appears in work of Garsia and Haiman
\cite{GH} on the bigraded character $DH_n(X;q,t)$ of the diagonal
harmonics. More specifically, we have
\[
DH_n(X;q,1/q)q^{{n \choose 2}}[n+1]_q = e_n[(1+q+\ldots+q^n)X].
\]
Curiously, ribbon functions also appear in recent work of \cite{HHLRU}.
Perhaps a deeper relationship between exists.

\textbf{Other generating functions.} In \cite{KLLT}, Kirillov,
Lascoux, Leclerc and Thibon gave a number of generating functions
for domino functions which were subsequently generalised in
\cite{Lam}.  As a special case, we have the following product
expansion for $n = 2$:
\[
\sum_\ll \g^{(2)}_\ll(X;q) = \frac{\prod_i(1+qx_i)}
{\prod_i(1-x_i)\prod_i(1-q^2x_i^2)\prod_{i<j}(1-x_ix_j)\prod_{i<j}(1-q^2x_ix_j)}.
\]
Can this be generalised to other values of $n$?

\textbf{Other incarnations of $\Lambda$.} Often the Fock space
$\f$ is identified with $\Lambda$ via
\[
\br{\ll} \leftrightarrow s_\ll.
\]
This gives $\f$ the extra structure of an algebra.  In this
context, our map $\m$ can be considered to be an operator
from $\Lambda(q)$ to $\Lambda(q)$.  In the notation of
\cite{LLT}, $\m$ would be the adjoint $\phi_q$ of the
operator $p_n^q$ which sends $h_\a$ to $\v_\a\cdot 1$.  It is not clear
whether this point of view leads to more results.

Leclerc \cite{Lec} has studied another embedding $\iota: \Lambda \rightarrow \f$ given by
\[
p_\ll \mapsto B_{-\ll}\br{0}.
\]
Altering this slightly, we may define a $\C(q)$-linear embedding
$\iota_q: \Lambda(q) \rightarrow \f$ by
\[
\plqn(p_\ll) \mapsto B_{-\ll}\br{0}.
\]
By Theorem \ref{thm:map}, we see that the composition
\[
\m \circ \iota_q: \Lambda(q) \rightarrow \Lambda(q)
\]
is the identity.
Leclerc has connected $\iota$ with the Macdonald polynomials and it is likely that
our setup can be connected with many other aspects of symmetric function theory
in this way.

\part{Combinatorics}
\label{part:comb}

\section{Ribbon Insertion}
\label{sec:ribinsertion} In this section we put the ribbon Pieri
formula (Theorem \ref{thm:pieri}) and ribbon Cauchy identity
(Theorem \ref{thm:cauchy}) in the context of ribbon
Robinson-Schensted-Knuth (RSK) insertion, where both will be
proven completely for the case $n=2$. Note that we will discuss
everything in the special case of an empty $n$-core, though
everything can be generalised to the non-empty $n$-core case.

\subsection{Robinson-Schensted-Knuth for usual Young tableaux}
Recall that the usual Robinson-Schensted bijection gives a
bijection between permutations $w \in S_m$ and pairs of standard
Young tableaux (see \cite{EC2}):
\[
w \mapsto \brac{P(w),Q(w)}.
\]
The tableaux $P(w)$ is defined recursively by the insertion
algorithm.  The tableaux $(T \leftarrow i)$ is obtained from $T$
by placing $i$ in the leftmost square of the first row of $T$
where the placement is legal ($i$ is placed legally if it is
greater than all numbers preceding it). If this displaces an
element $j$ then the element $j$ is placed in the leftmost square
of $T$ where the placement is legal. This continues until no more
bumping occurs. Then the insertion tableaux $P(w)$ is defined as
\[
P(w) = ((\cdots ((\emptyset \leftarrow w_1) \leftarrow w_2) \cdots
) \leftarrow w_n).
\]
The recording tableaux $Q(w)$ is defined by requiring that
\[
sh(Q(w)|_i) = sh((\cdots ((\emptyset \leftarrow w_1) \leftarrow
w_2) \cdots ) \leftarrow w_i),
\]
where $Q(w)|_i$ is the subtableaux of $Q(w)$ obtained by erasing
all the squares with numbers greater than $i$.

This bijection generalises to a bijection between two line arrays
\[
w = \brac{ \begin{array}{ccc} i_1 & \cdots & i_m \\ j_1 & \cdots &
j_m
\end{array}}
\]
satisfying (a) $i_1 \leq i_2 \leq \cdots \leq i_m$ and (b) if $i_r
= i_{r+1}$ then $j_r \leq j_{r+1}$.  The weight of the top row
then becomes the weight of the tableaux $Q(w)$ while the weight of
the bottom row becomes the weight of the tableaux $P(w)$.  This
immediately leads to the Cauchy identity (\ref{eq:classCauchy}) as
the generating function (where a number $i$ on the top row has
weight $y_i$ and a number $j$ on the borrom row has weight $x_j$)
for such two line arrays is exactly
\[
\prod_{i,j} \frac{1}{1-x_iy_j}.
\]
To obtain the bijection $w \mapsto (P(w),Q(w))$ one needs an
important property of insertion.  Let
\begin{align*}
T' &= (T \leftarrow i) \\
T'' &= (T' \leftarrow j)
\end{align*}
denote the the result from two successive insertions, and let
$\gamma = sh(T'/T)$ and $\theta = sh(T''/T)$ be the two new
squares added to the shapes.  Then $\gamma$ lies to the left of
$\theta$ if and only if $i \leq j$.  This \emph{increasing
insertion} property guarantees that $Q(w)$ will be semistandard.
In fact it is this property that is crucial to a combinatorial
proof (see \cite[p. 341]{EC2}) of the Pieri rule:
\[
h_ks_\ll = \sum_\mu s_\mu.
\]
We may interpret $h_k$ as the generating function for a $k$-tuple
of increasing positive integers $(i_1 \leq i_2 \leq \cdots \leq
i_k)$, and $s_\ll$ as the weight generating function of tableaux
$T$ with shape $\ll$, as usual. Then a bijection from the left
hand side to the right hand side is obtained by associating to a
pair $((i_1,\cdots,i_k),T)$ the tableaux
\[
T' = ((\cdots (( T \leftarrow i_1) \leftarrow i_2) \cdots )
\leftarrow i_k).
\]
The increasing insertion property guarantees that $sh(T')/\ll$ is
indeed a horizontal strip.

\subsection{General ribbon insertion}
\label{sec:genrib} Following \cite{SW2}, we will call a three line
array $\ww$ a $n$-colored biword if it is of the form
\[
\ww = \brac{ \begin{array}{ccc} c_1 & \cdots & c_m \\ i_1 & \cdots & i_m \\
j_1 & \cdots & j_m
\end{array}}
\]
where the $c_i$ are the `colors' taking values in
$\set{0,\ldots,n-1}$ and the $i_k$, $j_k$ are positive integers.
We will insist each such colored biword to have a canonical
ordering in such a way that only the multiset of ordered triples
$\set{(c_k,i_k,j_k)}$ matters.  For example, we could choose the
lexicographic ordering so (a) $c_k \leq c_{k+1}$ and (b) if $c_k =
c_{k+1}$ then $i_k \leq i_{k+1}$ and (c) if $c_k=c_{k+1}$ and $i_k
= i_{k+1}$ then $j_k \leq j_{k+1}$.  Giving the weight
\[
w((c_k,i_k,j_k)) = q^{2c_k}y_{i_k}x_{j_k}
\]
to each triple, the weight generating function of colored biwords
becomes
\[
\prod_{i,j} \prod_{k=0}^{n-1} \frac{1}{1-x_iy_jq^k}.
\]
Note that when $n = 1$ we just recover the setup of the previous
subsection.  We are then led to the following observation.
\begin{observation}
A (ribbon RSK) bijection $\pi: \ww \mapsto (P_r(\ww),Q_r(\ww))$
between \emph{colored biwords} and pairs of ribbon tableaux of the
same shape and fixed $n$-core will prove the Cauchy formula of
Theorem \ref{thm:cauchy} if
\begin{itemize}
\item
The bijection $\pi$ is weight preserving.  Thus the weight of the
second line of $\ww$ is $w(Q_r(\ww))$ and the weight of the third line
of $\ww$ is $w(P_r(\ww))$.
\item
The bijection $\pi$ sends color to spin.  Thus \be{eq:stc}
2(c_1+c_2+ \ldots + c_m) = s(P_r(\ww))+s(Q_r(\ww)). \ee
\end{itemize}
\end{observation}

Suppose now that the bijection $\pi$ is defined recursively via
insertion of ribbons $(c,j)$ into a tableaux $T$:
\[
P_r(\ww) = ((\cdots ((\emptyset \leftarrow (c_1,j_1)) \leftarrow
(c_2,j_2)) \cdots ) \leftarrow (c_m,j_m)).
\]
Then the tableaux $T' = T \leftarrow (c,j)$ should satisfy (a) the
tableaux $T'$ has an extra ribbon labelled $j$, (b) $sh(T')/sh(T)$
is a ribbon, and (c) $s(T') + s(sh(T')/sh(T)) = s(T) + 2c$. Here
we will think of $(c,j)$ as a ribbon labelled $j$ with spin $c$.
Let
\begin{align*}
T' &=T \leftarrow (c,j) \\
T'' &= T' \leftarrow (c',j'). \\
\end{align*}
A \emph{ribbon increasing insertion} property is a property of the
form
\begin{quote}
The ribbon $sh(T')/sh(T)$ lies to the left of $sh(T'')/sh(T')$ if
and only if $(c,j) \leq (c',j')$.
\end{quote}
Here $<$ should be some total order on ribbons labelled $j$ with
spin $c$.  Fix a ribbon tableaux $T$.  Then we can construct a bijection between sets of ribbons $\set{(c_i,j_i)}$ and ribbon tableaux $T'$ whose shape differs from that of $T$ by a horizontal ribbon strip by
\[
T' = ((\cdots((T\leftarrow (c_1,j_1)) \leftarrow (c_2,j_2))
\cdots) \leftarrow (c_k,j_k)).
\]
The ribbons $(c_i,j_i)$ are inserted according to the order $<$ thus ensuring the resulting shape changes by a horizontal ribbon strip.  Thus:

\begin{observation} A \emph{ribbon increasing insertion} property
for $\pi$ leads to a combinatorial proof of the ribbon Pieri formula
(Theorem \ref{thm:pieri}).
\end{observation}
Thus the generating function $H(t)$ of Section \ref{sec:pieri} can be
interpreted as the generating functions of ribbons $(c,j)$ with
weight $w(c,j) = q^{2k}x_j$.

\subsection{Domino insertion}
The above comments become proofs for the case $n=2$. Barbasch and
Vogan \cite{BV} have defined domino insertion in connection with
the primitive ideals of classical lie algebras.  This was put into
the usual bumping description by Garfinkle \cite{Gar}. Recently,
Shimozono and White \cite{SW} have extended Garfinkle's
description to the semistandard case and connected it with mixed
insertion.  They also observed that it had the crucial
color-to-spin property.  A straightforward extension to the
non-empty 2-core case was presented in \cite{Lam}.  We thus have:

\begin{thm}
Fix a 2-core $\d$. There is a bijection between colored biwords
$\ww$ of length $m$ with two colors $\set{0,1}$ and pairs
$(P_d(\ww),Q_d(\ww))$ of semistandard domino tableaux with the same
shape $\ll \in \p_{\d}$ and $|\ll| = 2m+|\d|$ with the following
properties:
\begin{enumerate}
\item
The bijection has the color-to-spin property:
\[
tc(\ww) = s(P_d(\ww)) + s(Q_d(\ww))
\]
where $tc(\ww)$ is the twice the sum of the colors in the top line
of $\ww$.
\item
The weight of $P_d(\ww)$ is the weight of the lowest line of $\ww$.
The weight of $Q_d(\ww)$ is the weight of the middle line of $\ww$.
\end{enumerate}
\end{thm}

In the standard case, Garfinkle's domino insertion is determined
by insisting that horizontal dominoes bump by rows and vertical
dominoes bump by columns.  More precisely, let $S$ be a domino
tableaux with no value repeated (but still semistandard), and $i$
some number not used in $S$.  We will describe $(S \leftarrow
(0,i))$ and $(S \leftarrow (1,i))$ which correspond to the
insertion of a horizontal (color 0) and vertical domino (color 1)
labelled $i$ respectively.

Let $T_{<i}$ be the subtableaux of $S$ consisting of all dominoes
labelled with numbers less than $i$.  Then set $T_{\leq i}$ to be
$T_{<i}$ union a horizontal domino in the first row labelled $i$
or a vertical domino in the first column labelled $i$ depending on
what we are inserting.  Now for $j >i$ we will recursively define
$T_{\leq j}$ given $T_{\leq j-1}$.  If there is no domino labelled
$j$ in $T$ then $T_{\leq j} = T_{\leq j-1}$.  Otherwise let
$\gamma_j$ denote the domino labelled $j$ in $S$ and set $\ll =
\text{sh}(T_{\leq j-1})$. We distinguish four cases.
\begin{enumerate}
\item
If $\gamma_j \cap \ll = \emptyset$ then set $T_{\leq j} = T_{\leq
j-1} \cup \gamma_j$.
\item
If $\gamma_j \cap \ll = \gamma_j$ is a horizontal domino in row
$k$ then $T_{\leq j}$ is obtained from $T_{\leq j-1}$ by adding a
horizontal domino labelled $j$ to row $k+1$.
\item
If $\gamma_j \cap \ll = \gamma_j$ is a vertical domino in column
$k$ then $T_{\leq j}$ is obtained from $T_{\leq j-1}$ by adding a
vertical domino labelled $j$ to column $k+1$.
\item
If $\gamma_j \cap \ll = (l,m)$ is a single square then $T_{\leq
j}$ is obtained from $T_{\leq j-1}$ by adding a domino labelled
$j$ so that the total shape of $T_{\leq j}$ is $\ll \cup (l+1,m+1)$.
\end{enumerate}
The resulting tableaux $T_{< \infty} = (S \leftarrow (c,i))$.

\begin{figure}[ht]
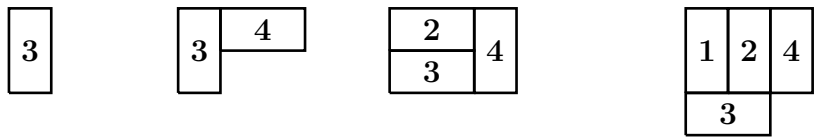

\pspicture(0,20)(304,80)
 \vdom(0,3){3}

\vdom(4,3){3} \hdom(5,4){4}

\hdom(9,4){2} \hdom(9,3){3} \vdom (11,3){4}

\vdom(16,3){1} \vdom(17,3){2} \hdom (16,2){3} \vdom(18,3){4}

\endpspicture
\caption{The result of the insertion $((((\emptyset \leftarrow (1,3))\leftarrow (0,4))\leftarrow (0,2) )\leftarrow (1,1))$.}
\label{fig:insertion}
\end{figure}

Figure \ref{fig:insertion} gives an example of domino insertion.

Shimozono and White's semistandard extension of domino insertion
leads automatically to the domino Cauchy formula as observed in
\cite{Lam}.  In \cite{Lam}, we have also described two dual domino
insertion algorithms which are bijections between `dual colored
biwords' and pairs of semistandard tableaux of conjugate shape.
This proves the dual domino Cauchy formula ($n=2$ in Theorem
\ref{prop:dualcauchy} here).

\medskip

It further turns out that Garfinkle's domino insertion has the
following domino increasing insertion property.  This was first
shown by Shimozono and White by connecting domino insertion with
mixed insertion. \cite{Lam} gives a different proof using growth
diagrams.  This domino increasing insertion property can be
described by specifying an order $<$ on dominoes as follows
($\gamma_i$ denotes a domino labelled $i$)
\begin{enumerate}
\item
If $\gamma_i$ is horizontal and $\gamma_j$ vertical then
$\gamma_i
> \gamma_j$.
\item
If $\gamma_i$ and $\gamma_j$ are both horizontal then $\gamma_i >
\gamma_j$ if and only if $i > j$.
\item
If $\gamma_i$ and $\gamma_j$ are both vertical then $\gamma_i >
\gamma_j$ if and only if $i < j$.
\end{enumerate}
Under this order, Garfinkle's domino insertion has a ribbon
increasing insertion property, as described in Section
\ref{sec:genrib}:
\begin{lem}
Let $T$ be a domino tableaux without the labels $i$ and $j$.  Set
$T' = (T \leftarrow \gamma_i)$ and $T'' = (T' \leftarrow
\gamma_j)$ for some dominoes $\gamma_i$ and $\gamma_j$.  Then
$sh(T'/T)$ lies to the left of $sh(T''/T')$ if and only if
$\gamma_i < \gamma_j$.
\end{lem}

Similarly, the dual domino insertion has a property which is dual
to this.  This increasing property is retained when the bijection
is extended to the semistandard case which we shall now describe
in brief (see \cite{SW,Lam} for details). Let $T$ be a
semistandard domino tableaux and $(c,j)$ a domino we want to
insert, where $j$ is a value possibly occurring in $T$. For each
value $i$, all the dominoes labelled $i$ in $T$ can be ordered
from left to right and labelled $i_1, i_2, \ldots$, where $i-1<
i_1 < i_2 < \cdots < i+1$ -- all $i_a$ behave like an `$i$' when
compared to any other value. Then this insertion can be simulated
by treating the new label $j$ as being larger than or smaller than
all other values $j_a$ (present in $T$) depending on whether $c
=0$ or $c=1$.  We may then perform insertion as in the standard
case.  Afterwards we rename all the $i_a$ to $i$ to obtain $(T
\leftarrow (c,j))$.  In particular, one may check that the
increasing insertion property is compatible with the semistandard
insertion.

Immediately we obtain
\begin{prop}
Semistandard domino insertion gives a combinatorial proof of the
Pieri rule (Theorem \ref{thm:pieri}) for $n=2$.  Dual semistandard
domino insertion gives a combinatorial proof of the dual Pieri
rule for $n=2$.
\end{prop}

\subsection{Shimozono and White's ribbon insertion}
Shimozono and White \cite{SW2} have described a ribbon insertion
algorithm for general $n$.  This can be described in a traditional bumping fashion
or in terms of Fomin's growth diagrams \cite{Fom1,Fom2}.

The ribbon insertion algorithm of \cite{SW2} has the usual weight
preserving properties, but also the spin to color property
(\ref{eq:stc}) which an earlier ribbon-RSK algorithm of Stanton
and White \cite{SW1} did not have.  However, the algorithm stops
short of being a bijection between colored biwords and pairs of
semistandard ribbon tableaux.  The algorithm is only described as
a bijection $\pi$ between colored words $\ww$ (not biwords) and a pair
$(P_r(\ww),Q_r(\ww))$ where $P_r(\ww)$ is a semistandard ribbon tableaux and
$Q_r(\ww)$ is a standard ribbon tableaux.  In particular the Cauchy
identity of Theorem \ref{thm:cauchy} does not immediately follow.
The algorithm also does not seem to possess a \emph{ribbon
increasing insertion} property.  However one can at least salvage
the following, which is just the first Pieri rule.

\begin{prop}
Shimozono and White's bijection $\pi$ gives a combinatorial proof
that
\[
(1+q^2 +\ldots + q^{2(n-1)})h_1 \g_\ll(X;q) = \sum_\mu
q^{s(\mu/\ll)}\g_\mu(X;q)
\]
where the sum is over all $\mu$ such that $\mu/\ll$ is a
$n$-ribbon.
\end{prop}
\begin{proof}
As before we construct a weight preserving bijection between the two
sides of the Pieri rule by:
\[
(T,(c,j)) \mapsto T' = (T \leftarrow (c,j)).
\]
The color $c$ ranges from 0 to $n-1$ and $h_1$ is just the generating function
for the labels $j$.
\end{proof}

Shimozono and White's ribbon insertion is determined by forcing
all ribbons to bump by rows to another ribbon of the same spin (at
least in the standard case).  It is possible however to insist
that all ribbons of a particular spin bump by columns instead.
Unfortunately, it appears that none of these algorithms have a
ribbon increasing insertion property.

\section{Murnagham-Nakayama and Pieri}
\label{sec:mnviapieri}
In this (self-contained) section we will study the formal
combinatorial relationship between Murnagham-Nakayama and Pieri
rules for ribbon tableaux.

In particular we will obtain a direct proof that the usual
Murnagham-Nakayama rule and Pieri rules are formally equivalent in
a combinatorial fashion. This bypasses the usual method of proof
which goes via the Jacobi-Trudi formulae.  The only algebraic fact
needed is the following lemma:

\begin{lem}
\label{lem:phe}
The power sum and elementary symmetric functions satisfy the
following equation
\[ ne_n = p_1e_{n-1} - p_2e_{n-2} + \cdots + (-1)^{n-1} p_n.\]  Similarly,
we have
\[ mh_m = p_{m-1}h_1 + p_{m-2}h_2 + \cdots + p_m\] for the homogenoeous
and power sum symmetric functions.
\end{lem}
\begin{proof}
See (2.10) in \cite{Mac}.
\end{proof}

Let $V$ be a vector space over $\C(q)$ and $v_\ll$ be vectors in
$V$ labelled by partitions.  Recall the definitions of $\xx^{\mu/\ll}_k(q)$,
$\kk_{\mu/\ll,k}(q)$ and $\l_{\mu/\ll,k}(q)$ from Section \ref{sec:not}.
Suppose $\set{P_k}$ are commuting
linear operators satisfying
\[
P_k v_\ll = \sum_\mu \xx_k^{\mu/\ll}(q) v_\mu \;\;\; \mbox{for all
$k$}
\]
then we will say that the Murnagham-Nakayama rule holds.
\par
Suppose $\set{H_k}$ are commuting linear operators on $V$ satisfying
\[
H_k v_\ll = \sum_\mu \kk_{\mu/\ll,k}(q) v_\mu \;\;\; \mbox{for all
$k$},
\]
then we will say that Pieri formula holds.
\par
Suppose $\set{E_k}$ are commuting linear operators on $V$ satisfying
\[
E_k v_\ll = \sum_\mu \l_{\mu/\ll,k}(q) v_\mu \;\;\; \mbox{for all
$k$},
\]
then we will say that dual-Pieri formula holds.

If the skew shapes $\mu/\ll$ are replaced by $\ll/\mu$ in the
above formulae, we get adjoint versions of these formulae which
can be thought of as lowering operators.  Thus if a set of
commuting linear operators $\set{P_k^{\perp}}$ satisfies
\[
P_k^\perp v_\ll = \sum_\mu \xx_k^{\ll/\mu}(q) v_\mu \;\;\;
\mbox{for all $k$}
\]
then we will say the lowering Murnagham-Nakayama rule holds, and
similarly for $\set{E_k^\perp}$ and $\set{H_k^\perp}$.

\begin{prop}
\label{prop:mnviapieri1} Fix $n \geq 1$ as usual.
Let $\set{H_k}$ and $\set{P_k}$ be commuting sets of linear operators
satisfying the relations between $h_k$ and $p_k$ in $\Lambda$.  Then
the ribbon Murnagham-Nakayama rule holds for $\set{P_k}$ if and only if
the ribbon Pieri rule holds for $\set{H_k}$.
\end{prop}

\begin{proof}
Let us suppose the Murnagham-Nakayama rule holds for $\set{P_k}$.
We will proceed by induction on $k$.  Since $H_1 = P_1$ the
starting condition is clear.  Now suppose the proposition has been
shown up to $k-1$.  By Lemma \ref{lem:phe}, $kH_k$ acts on $v_\ll$
in the same way that $H_{k-1}P_1 + H_{k-2}P_2 + \cdots + P_k$
does.

The action of the latter on $v_\ll$ gives some linear combination
of $v_\mu$ where $\mu/\ll$ is a union $(S,T)$ of a ribbon border
strip $S$ and a horizontal ribbon strip $T$ (where $T$ is a
horizontal ribbon strip of the shape $\mu/(\ll \cup S)$).  Denote
by $S_1, \ldots, S_a$ the distinguished decomposition of $S$ into
horizontal ribbon strips.  Fixing $\mu$ we now consider the set
$\S$ of such ordered pairs $(S,T)$ where $S$ has size between $1$
and $k$, while $T$ has size from $k-1$ to $0$.

We now place the $(S,T)$ into equivalence classes $\P(S,T)$.  The
equivalence relation is given by taking the transitive closure of
the relation \begin{equation} \label{eq:eq} (S,T)\sim (S-S_a,T
\cup S_a)\end{equation} for every pair such that $T \cup S_a$ is a
horizontal strip. This relation is ill-defined when $a = 1$, that
is when $S$ is actually a single connected horizontal ribbon
strip (in which case $S-S_1$ is empty and $(S-S_1,T \cup S_1)$
does not belong to $\S$), which we shall ignore for the moment.

Let us consider any other equivalence class $\P = \P(S,T)$.  We
claim that it contains a unique element $(S',T')$ (where $S' =
\set{S_1',\ldots,S_a'}$) such that $T' \cup S_a'$ is not a
horizontal ribbon strip.  This is due to the definition of a
border ribbon strip which ensures that the right hand side of
(\ref{eq:eq}) always has this property.  Thus the graph of the
relations (\ref{eq:eq}) is star-shaped, proving our claim.  Now
let $C$ be a component of $T'$ such that $C \cup S_a'$ is not a
horizontal ribbon strip. Then there is a unique sub-horizontal
ribbon strip $C'$ of $C$ which can be added to $S'$ to form a
ribbon strip. This $C'$ may be described as follows. Order the
ribbons of $C$ from left to right $c_1, c_2, \ldots, c_l$.  Find
the smallest $i$ such that $c_i$ touches the bottom of $S_a'$ and
we set $C' = \set{c_1,c_2, \ldots, c_i}$.  We call $C$ a critical
component and $C'$ the nice part of $C$.

Then the equivalence class $\P$ is exactly $(S',T')$ together with
the pairs $(S,T)$ such that $S = \set{S_1,\ldots,S_{a+1}}$ where
$S_i = S_i'$ for $1 \leq i \leq a$, and $S_{a+1}$ is the union of
the nice parts of some (arbitrary) subset of the set critical
components of $T'$.  It is immediate from the construction that
$(S,T)$ will be a valid pair in $\S$.  We observe that the
contribution of $\P$ 
\[
\sum_{(S,T) \in \P} (-1)^{h(S)} q^{s(S \cup T)}
\]
to the coefficient of $v_\mu$ is exactly 0,
since the the tiling and hence the spin of the contribution is
fixed and the definition of height is exactly so that the signs
sum up to 0 (this corresponds to the identity $(1-1)^c = 0$).

It remains to consider the elements $(S,T)$ where $S$ is a
connected horizontal ribbon strip such that $S \cup T$ is also a
horizontal ribbon strip.  Since $S$ is connected we can recover it
from $S \cup T$ by specifying its rightmost ribbon.  Thus such
pairs occur exactly $k$ times for each horizontal ribbon strip of
shape $\mu/\ll$, and hence the Pieri rule is satisfied for the
operator $H_k$.

The converse clearly follows from the same argument.
\end{proof}

\medskip
\begin{thm}
\label{thm:mnviapieri} Let $\set{H_i}$, $\set{E_i}$ and
$\set{P_i}$ be commuting operators on a vector space $V$ over
$\C(q)$ satisfy the relations of $h_i$, $e_i$ and $p_i$ in
$\Lambda$.  Let $v_\ll$ be a set of vectors in $V$ indexed by
partitions.  Suppose that one of the Pieri, dual-Pieri and
Murnagham-Nakayama holds, then all three holds.  The same is true
for the lowering operators satisfying the same relation.
\end{thm}
\begin{proof}
That the Murnagham-Nakayama rule and Pieri rules are equivalent is
just Proposition \ref{prop:mnviapieri1}. One way to see that the
Pieri rules and dual-Pieri rules are equivalent is to argue as before
and use the
relation
\[
h_m - h_{m-1}e_1 + \cdots + (-1)^m e_m = 0
\]
which is easily deduced from the generating functions $H(t) =
\sum_m h_m t^m$ and $E(t) = \sum_m e_m t^m$. However, a short cut
is to use Proposition \ref{prop:hor}.  We see that both the Pieri
and dual-Pieri formulae hold in $\f$ for operators satisfying the
relations of $h_i$ and $e_i$.  But the vectors $\br{\ll}$ are
linearly independent in $\f$ so this formally implies (by
linearity) that the same is true for any set of vectors $v_\ll$ in
a vector space $V$ over $\C(q)$.

The argument is identical for lowering operators.
\end{proof}

Note that the condition on a horizontal ribbon strip to be
connected can be described in terms of the $n$-quotient as
follows. Let $T$ be a ribbon tableaux with $n$-quotient
$\set{T^{(0)}, \ldots,T^{(n-1)}}$. Let $\set{(d_i, p_i)}$ be the
set of diagonals which are nonempty in the $n$-quotient of the
horizontal ribbon strip $R$. Thus diagonal $diag_{d_i}$ of
$T^{(p_i)}$ contains a square corresponding to some ribbon in the
horizontal ribbon strip $R$.  Then the horizontal ribbon strip $R$
is connected if and only if the set of integers $\set{d_i}$ is an
interval (connected) in $\zz$.  Thus border ribbon strips may be
characterised in terms of the $n$-quotient.

\end{document}